\documentclass[11pt]{amsart}
\usepackage{amssymb}
\usepackage{amscd}		% commutative diagram

   \title[Vanishing theorems and character formulas]{Vanishing
          theorems and character formulas for the Hilbert scheme of
          points in the plane}
                                                                   
   \author{Mark Haiman}                                            
                                                                   
   \address{Dept.\ of Mathematics\\
            University of California\\
            Berkeley, CA, 94720-3840}

   \email{mhaiman@math.berkeley.edu}

   \date{Jan 14, 2002}                                        

\thanks{Research supported in part by N.S.F. grants DMS-9701218 and
DMS-0070772 and the Isaac Newton Institute}

\DeclareMathOperator{\Hilb}{Hilb}
\DeclareMathOperator{\Spec}{Spec}
\DeclareMathOperator{\Hom}{Hom}
\DeclareMathOperator{\Supp}{Supp}
\DeclareMathOperator{\GL}{GL}
\DeclareMathOperator{\SL}{SL}
\DeclareMathOperator{\ch}{char}
\DeclareMathOperator{\tr}{tr}

\DeclareMathOperator{\ev}{ev}
\DeclareMathOperator{\Dh}{DH}
\DeclareMathOperator{\PF}{PF}
\DeclareMathOperator{\OP}{OP}

\newcommand{\Ocal}{{\mathcal O}}
\newcommand{\Fcal}{{\mathcal F}}
\newcommand{\Hcal}{{\mathcal H}}

\newcommand{\Bcal}{{\mathcal B}}

\newcommand{\Htild}{\tilde{H}}
\newcommand{\Ktild}{\tilde{K}}

\newcommand{\mfrak}{{\mathfrak m}}

\newcommand{\ZZ}{{\mathbb Z}}
\newcommand{\QQ}{{\mathbb Q}}
\newcommand{\NN}{{\mathbb N}}
\newcommand{\CC}{{\mathbb C}}
\newcommand{\TT}{{\mathbb T}}

\renewcommand{\aa}{{\mathbf a}}
\newcommand{\bb}{{\mathbf b}}
\newcommand{\xx}{{\mathbf x}}
\newcommand{\yy}{{\mathbf y}}
\newcommand{\zz}{{\mathbf z}}
\newcommand{\ee}{{\mathbf e}}
\newcommand{\boldt}{{\mathbf t}}
\newcommand{\ggam}{\text{\boldmath $\gamma $}}

\newcommand{\defeq}{\;\underset{\operatorname{def}}{=}\;}

\newcommand{\smallfrac}[2]{{\textstyle \frac{#1}{#2}}}

\newcommand{\ext}{\wedge}

\newcommand{\subdot}{\text{{\bf .}}}

\newcommand{\llbracket}{{[\![}}
\newcommand{\rrbracket}{{]\!]}}
\newcommand{\sslash}{{/\!/}}

\theoremstyle{plain}
\newtheorem{prop}{Proposition}[section]
\newtheorem{thm}[prop]{Theorem}
\newtheorem{cor}[prop]{Corollary}
\newtheorem{lem}[prop]{Lemma}
\newtheorem{conj}[prop]{Conjecture}

\theoremstyle{definition}

\begin{document}

\begin{abstract}
In an earlier paper \cite{Hai01}, we showed that the Hilbert scheme of
points in the plane $H_{n} = \Hilb ^{n}(\CC ^{2})$ can be identified
with the Hilbert scheme of regular orbits $\CC ^{2n} \sslash S_{n}$.
Using this result, together with a recent theorem of Bridgeland, King
and Reid \cite{BrKiRe01} on the generalized McKay correspondence, we
prove vanishing theorems for tensor powers of tautological bundles on
the Hilbert scheme.  We apply the vanishing theorems to establish
(among other things) the character formula for diagonal harmonics
conjectured by Garsia and the author in \cite{GaHa96}.  In particular
we prove that the dimension of the space of diagonal harmonics is
equal to $(n+1)^{n-1}$.
\end{abstract}

\maketitle


\section{Introduction}
\label{s:intro}

In this article we continue the investigation begun in \cite{Hai01} of
the geometry of the Hilbert scheme of points in the plane and its
algebraic and combinatorial implications.  In the earlier article, we
showed that the isospectral Hilbert scheme has Gorenstein
singularities, thereby proving the ``$n!$ conjecture'' of Garsia and
the author and the positivity conjecture for Macdonald polynomials.
Here we extend these results by proving vanishing theorems for tensor
products of tautological vector bundles over the Hilbert scheme $H_{n}
= \Hilb ^{n}(\CC ^{2})$ and its zero fiber $Z_{n}$ (the preimage of
the origin under the Chow morphism $\sigma \colon H_{n}\rightarrow
S^{n}\CC ^{2}$).

The algebraic-combinatorial consequence of the new results is a series
of character formulas for the spaces of global sections of the vector
bundles in question.  As a particular case, we obtain the character
formula for the ring of coinvariants of the diagonal action of the
symmetric group $S_{n}$ on $\CC ^{2n}$, or equivalently, for the space
of {\it diagonal harmonics}.  This character formula had been
conjectured by Garsia and the author in \cite{GaHa96}, where we showed
that it implies several earlier conjectures in \cite{Hai94} relating
the character of the diagonal harmonics to $q$-Lagrange inversion,
$q$-Catalan numbers, and $q$-enumeration of rooted forests and parking
functions.  One of these earlier conjectures, now proven, is that the
dimension of the space of diagonal harmonics is
\begin{equation}\label{e:DHn-formula}
\dim \Dh _{n} = (n+1)^{n-1}.
\end{equation}
Another is that the Hilbert series of the doubly-graded space $(\Dh
_{n})^{\epsilon }$ of $S_{n}$-alternating diagonal harmonics is given
by the $q,t$-Catalan polynomial $C_{n}(q,t)$ from \cite{GaHa96,Hai98}.
Hence $C_{n}(q,t)$ has positive integer coefficients.  Recently,
Garsia and Haglund \cite{GaHag01} gave a different proof of this fact,
based on a combinatorial interpretation of the coefficients.  Yet
another conjecture in \cite{Hai94} was that the space of diagonal
harmonics is generated by certain $S_{n}$-invariant polarization
operators applied to the space of classical harmonics.  We prove this
``operator conjecture'' here, using our identification of the
coinvariant ring with the space of global sections of a vector bundle
on $Z_{n}$.


To describe our results further, we recall from \cite{Hai01} that
$H_{n}$ is isomorphic to the Hilbert scheme of orbits $\CC ^{2n}
\sslash S_{n}$ for the diagonal action of $S_{n}$ on $\CC ^{2n}$.
Full definitions will be given in Section~\ref{s:main}; for now we
merely fix terminology in order to announce our main theorems.  On the
Hilbert scheme $H_{n}$ we have a natural tautological vector bundle
$B$ of rank $n$, while on $\CC ^{2n} \sslash S_{n}$ we have a
tautological bundle $P$ of rank $n!$, with an $S_{n}$ action in which
each fiber affords the regular representation.  We can view both $B$
and $P$ as bundles on $H_{n}$ via the isomorphism $H_{n}\cong \CC
^{2n} \sslash S_{n}$.  The usual tautological bundle $B$ is the
pushdown to $H_{n}$ of the sheaf $\Ocal _{F_{n}}$ of regular functions
on the universal family $F_{n}$ over $H_{n}$.  The ``unusual''
tautological bundle $P$ may similarly be identified with the pushdown
of the sheaf $\Ocal _{X_{n}}$ of regular functions on the isospectral
Hilbert scheme $X_{n}$, which is actually the universal family over
$\CC ^{2n}\sslash S_{n}$.

Our first main result, Theorem~\ref{t:van}, is a vanishing theorem for
the higher cohomology groups $H^{i}(H_{n},{P\otimes B^{\otimes l}})$,
$i>0$ of the tensor product of $P$ with any tensor power of $B$.  We
also identify the space of global sections $H^{0}(H_{n},{P\otimes
B^{\otimes l}})$.  The latter turns out to be the coordinate ring
$R(n,l)$ of the polygraph, a subspace arrangement defined
in \cite{Hai01}, which plays an important technical role there and
again here.  This identification of $R(n,l)$ with
$H^{0}(H_{n},{P\otimes B^{\otimes l}})$ explains why the polygraph
carries geometric information about the Hilbert scheme, an explanation
which we were only able to hint at in \cite{Hai01}.  Our theorem
extends vanishing theorems of Danila \cite{Dan01} for the tautological
bundle $B$ and of Kumar and Thomsen \cite{KuTh01} for the natural
ample line bundles $\Ocal_{H_{n}} (k)$, $k>0$.  Indeed, it implies the
vanishing of the higher cohomology groups $H^{i}(H_{n}, {\Ocal
(k)\otimes B^{\otimes l}})$ for all $k,l\geq 0$.  This is an immediate
corollary, since the trivial bundle $\Ocal _{H_{n}}$ is a direct
summand of $P$, and the line bundle $\Ocal _{H_{n}}(1)$ is the highest
exterior power of $B$.


Our second main result, Theorem~\ref{t:van-zero}, is a vanishing
theorem for the same vector bundles on the zero fiber $Z_{n}$.  The
vanishing part of this second theorem follows immediately from the
first theorem, applied to an explicit locally free resolution of
$\Ocal _{Z_{n}}$ described in \cite{Hai98} and reviewed in detail in
Section~\ref{s:main}, below.  By examining the resolution more
closely, we can also identify the space of global sections
$H^{0}(Z_{n},{P\otimes B^{\otimes l}})$.  When $l=0$ it turns out that
$H^{0}(Z_{n},P)$ coincides with the coinvariant ring for the diagonal
$S_{n}$ action on $\CC ^{2n}$, yielding the applications to diagonal
harmonics.

Character formulas for the spaces of global sections follow from our
vanishing theorems by an application of the Atiyah--Bott Lefschetz
formula \cite{AtBo66}.  For the diagonal harmonics, the calculation
completes a program outlined by Procesi, who first proposed this
method of determining the character.  To carry out the calculation, we
need to know the characters of the fibers of $P$ at distinguished
torus-fixed points $I_{\mu }$ on $H_{n}$.  By our results in
\cite{Hai01}, these characters are given by the Macdonald polynomials.
The character formulas we obtain here are therefore expressed in terms
of Macdonald polynomials.  Specifically, they are symmetric functions
with coefficients depending on two parameters $q, t$.  By virtue of
being characters, these symmetric functions are necessarily {\it
$q,t$-Schur positive}, that is, they are linear combinations of Schur
functions by polynomials or power series in $q$ and $t$ with positive
integer coefficients.  It develops that certain operator expressions
considered in \cite{BeGaHaTe99} are instances of our character
formulas, whose positivity partially establishes \cite[Conjecture
V]{BeGaHaTe99}.  The full conjecture in \cite{BeGaHaTe99} is slightly
stronger than what we obtain here.  Its proof using the methods of
this paper would require an improved vanishing theorem, which we offer
as a conjecture at the end of Section~\ref{s:char}.


Among our character formulas is one for the polygraph coordinate ring
$R(n,l)$ as a doubly graded algebra.  Specializing this, we get a
formula for its Hilbert series $\Hcal _{R(n,l)}(q,t)$ in terms of
symmetric function operators whose eigenfunctions are Macdonald
polynomials.  A combinatorial interpretation of $\Hcal _{R(n,l)}(q,t)$
is implicit in the basis construction for $R(n,l)$ in \cite{Hai01}.
It can be made explicit (although we will not do so here), yielding an
identity between a combinatorial generating function and the
expression involving Macdonald operators in Corollary~\ref{c:HR(n,l)},
below.  This is one of only two combinatorial interpretations known at
present for $q,t$-(Schur) positive expressions arising from our
character formulas.  The other is the Garsia--Haglund interpretation
of $C_{n}(q,t)$ mentioned above.  An important problem that remains
open is to combinatorialize all the character formulas.


In Section~\ref{s:main}, we give definitions and state our two main
theorems in full.  We then apply Theorem~\ref{t:van} to deduce
Theorem~\ref{t:van-zero}.  We deduce the character formulas and the
operator conjecture from the vanishing theorems in
Sections~\ref{s:char} and \ref{s:op}.  For the proof of
Theorem~\ref{t:van}, we combine results from \cite{Hai01} with a
recent general theorem of Bridgeland, King and Reid \cite{BrKiRe01}.
This is done in Section~\ref{s:van}.  To complete this introduction,
we preview the proof of Theorem~\ref{t:van}.

The Bridgeland--King--Reid theorem concerns the Hilbert scheme of
orbits $V\sslash G$, for a finite subgroup $G\subseteq \SL (V)$.  The
theorem has two parts.  The first part (which we will not use) is a
criterion for $V\sslash G$ to be a {\it crepant resolution of
singularities} of $V/G$, meaning that $V\sslash G$ is non-singular and
its canonical sheaf is trivial.  The second (and for us, crucial) part
says that when the criterion holds there is an equivalence of
categories $\Phi \colon D(V\sslash G)\rightarrow D^{G}(V)$.  Here
$D(V\sslash G)$ is the derived category of complexes of sheaves on
$V\sslash G$ with bounded, coherent cohomology, and $D^{G}(V)$ is the
similar derived category of $G$-equivariant sheaves on $V$.


Our identification of $\CC ^{2n}\sslash S_{n}$ with $H_{n}$ shows that
the Bridgeland--King--Reid criterion holds for $V=\CC ^{2n}$,
$G=S_{n}$.  It is well-known that $H_{n}$ is a crepant resolution of
$\CC ^{2n}/S_{n} = S^{n}\CC ^{2}$, which is why we don't need the
first part of their theorem.  By the second part, however, we have an
equivalence $\Phi $ between the derived category $D(H_{n})$ of sheaves
on the Hilbert scheme and the derived category $D^{S_{n}}(\CC ^{2n})$
of finitely generated $S_{n}$-equivariant modules over the polynomial
ring $\CC [\xx ,\yy ]$ in $2n$ variables.  In this notation,
Theorem~\ref{t:van} reduces to an identity $\Phi B^{\otimes l} =
R(n,l)$.  Denoting the inverse equivalence by $\Psi $, we may rewrite
this as $\Psi R(n,l) = B^{\otimes l}$, which is the form in which we
prove it.  The advantage of this form is that there is no sheaf
cohomology involved in the calculation of $\Psi $, only commutative
algebra.  Conveniently, the commutative algebraic fact we need is
precisely the freeness theorem for the polygraph ring $R(n,l)$, which
was the key technical theorem in \cite{Hai01}.  Thus we use here both
the geometric results from \cite{Hai01} and the main algebraic
ingredient in their proof.


In closing, let us remark that a number of important problems relating
to this circle of ideas remain open.  We have already mentioned the
problem of combinatorializing the rest of the character formulas.
Another set of problems involves phenomena in three or more sets of
variables.  We expect, for example, that the analog of the operator
conjecture should continue to hold in additional sets of variables
$\xx ,\yy ,\ldots,\zz $.  For exactly three sets of variables, we
remind the reader of the empirical conjecture in \cite{Hai94} that the
dimension of the space of ``triagonal'' harmonics should be
\begin{equation}\label{e:THn-formula}
2^{n}(n+1)^{n-2},
\end{equation}
and that of its $S_{n}$-alternating subspace should be
\begin{equation}\label{e:tri-Catalan}
(3n+3)(3n+4)\cdots (4n+1)/3\cdot 4\cdots (n+1).
\end{equation}
Our present methods do not readily apply to these problems, as we make
heavy use of special properties of the Hilbert scheme $\Hilb ^{n}(\CC
^{2})$ that do not hold for $\Hilb ^{n}(\CC ^{d})$ with $d\geq 3$.
Another open problem is to generalize from $S_{n}$ to other Weyl
groups or complex reflection groups.  Such a generalization will not
be entirely straightforward, as shown by some obstacles discussed in
\cite{Hai94} and \cite{Hai01}.  Finally, despite the strength of the
vanishing theorems proven here, they surely are not the strongest
possible.  The conjecture at the end of Section~\ref{s:char} suggests
one possible improvement.


\section{Definitions and main theorems}
\label{s:main}

We denote by $H_{n}$ the Hilbert scheme of points $\Hilb ^{n}(\CC
^{2})$ parametrizing $0$-dimensional subschemes of length $n$ in the
affine plane over $\CC $.  By Fogarty's theorem \cite{Fog68}, $H_{n}$
is irreducible and non-singular, of dimension $2n$.  As a matter of
notation, if $V(I)\subseteq \CC ^{2}$ is the subscheme corresponding
to a (closed) point of $H_{n}$, we refer to this point by its defining
ideal $I\subseteq \CC [x,y]$.  Thus $H_{n}$ is identified with the set
of ideals $I$ such that $\CC [x,y]/I$ has dimension $n$ as a complex
vector space.

The {\it multiplicity} of a point $P\in V(I)$ is the length of the
Artin local ring $(\CC [x,y]/I)_{P}$.  The multiplicities of all
points in $V(I)$ sum to $n$, giving rise to a $0$-dimensional
algebraic cycle $\sum _{i}m_{i}P_{i}$ of weight $\sum _{i}m_{i} = n$.
We may view this cycle as an unordered $n$-tuple $\llbracket
P_{1},\ldots,P_{n} \rrbracket \in S^{n}\CC ^{2}$, in which each point
is repeated according to its multiplicity.  The {\it Chow morphism}
\begin{equation}\label{e:Chow-Hn}
\sigma \colon H_{n}\rightarrow S^{n}\CC ^{2} = \CC ^{2n}/S_{n}
\end{equation}
is the projective and birational morphism mapping each $I\in H_{n}$ to
the corresponding algebraic cycle $\sigma (I) = \llbracket
P_{1},\ldots,P_{n} \rrbracket $.


We denote by $F_{n}$ the {\it universal family} over the Hilbert scheme,
\begin{equation}\label{e:F->Hn}
\begin{CD}
F_{n}&	\; \subseteq \; H_{n}\times \CC ^{2}\\
@V{\pi }VV\\
H_{n},
\end{CD}
\end{equation}
whose fiber over a point $I\in H_{n}$ is the subscheme $V(I)\subseteq
\CC ^{2}$.  The universal family is flat and finite of degree $n$ over
$H_{n}$, and hence is given by $F_{n} = \Spec B$, where $B = \pi
_{*}\Ocal _{F_{n}}$ is a locally free sheaf of $\Ocal
_{H_{n}}$-algebras of rank $n$.  Here and elsewhere we identify any
locally free sheaf of rank $r$ with the rank $r$ algebraic vector
bundle whose sheaf of sections it is.  Then $B$ is the {\it
tautological vector bundle}, the quotient of the trivial bundle ${\CC
[x,y] \otimes _{\CC } \Ocal_{H_{n}}}$ with fiber $\CC [x,y]/I$ at each
point $I\in H_{n}$.


If $G$ is a finite subgroup of $\GL (V)$, where $V = \CC ^{d}$ is a
finite-dimensional complex vector space, we denote by $V \sslash G$
the Hilbert scheme of regular $G$-orbits in $V$, as defined by Ito and
Nakamura \cite{ItNa96,ItNa99}.  Specifically, if $v\in V$ has trivial
stabilizer (as is true for all $v$ in a Zariski open set), then its
orbit $Gv$ is a point of $\Hilb ^{|G|}(V)$, and $V\sslash G$ is the
closure in $\Hilb ^{|G|}(V)$ of the locus of all such points.  By
definition, $V \sslash G$ is irreducible.  The universal family over
$\Hilb ^{|G|}(V)$ restricts to a universal family
\begin{equation}\label{e:X->V//G}
\begin{CD}
X&	\; \subseteq \; (V\sslash G) \times V\\
@V{\rho }VV\\
V\sslash G.
\end{CD}
\end{equation}
The group $G$ acts on $X$ and on the tautological bundle $P = \rho
_{*}\Ocal _{X}$.  This action makes $P$ a vector bundle of rank $|G|$
whose fibers afford the regular representation of $G$.  There is a
canonical Chow morphism $V\sslash G\rightarrow V/G$, which can be
conveniently defined as follows.  Since $P$ is a sheaf of $\Ocal
_{V\sslash G}$-algebras, it comes equipped with a homomorphism $\Ocal
_{V\sslash G} \rightarrow P$.  This homomorphism is an isomorphism of
$\Ocal _{V\sslash G}$ onto the sheaf of invariants $P^{G}$.
Geometrically, this means that the map $X/G\rightarrow V\sslash G$
induced by $\rho $ is an isomorphism.  The canonical projection
$X\rightarrow V$ induces a morphism $X/G\rightarrow V/G$ whose
composite with the isomorphism $V\sslash G\cong X/G$ yields the Chow
morphism.  The Chow morphism is projective and birational, restricting
to an isomorphism on the open locus consisting of orbits $Gv$ for $v$
with trivial stabilizer.


The case of interest to us is $V = \CC ^{2n}$, $G = S_{n}$, where
$S_{n}$ acts on $\CC ^{2n} = (\CC ^{2})^{n}$ by permuting the
cartesian factors.  This is the same as the diagonal action of $S_{n}$
on the direct sum of two copies of its natural representation $\CC
^{n}$.  Coordinates on $\CC ^{2n}$ will be denoted
\begin{equation}\label{e:C^2n-coords}
\xx , \yy  = x_{1},y_{1},\ldots,x_{n},y_{n};
\end{equation}
then $S_{n}$ acts by permuting the $x$ variables and the $y$ variables
simultaneously.  In \cite{Hai99} we constructed a canonical morphism
$\CC ^{2n}\sslash S_{n}\rightarrow H_{n}$ such that the composite
\begin{equation}\label{e:Chow-C2n//Sn}
\CC ^{2n}\sslash S_{n}\rightarrow H_{n}\overset{\sigma }{\rightarrow
} S^{n}\CC ^{2}
\end{equation}
is the Chow morphism for $\CC ^{2n}\sslash S_{n}$.  By Theorem~5.1
of \cite{Hai01}, the canonical morphism is an isomorphism $\CC
^{2n}\sslash S_{n}\cong H_{n}$.

The universal family over $\CC ^{2n}\sslash S_{n}$ will be denoted
$X_{n}$.  We identify $\CC ^{2n}\sslash S_{n}$ with $H_{n}$ by means
of the canonical isomorphism, so that the projection $\rho $ of the
universal family onto $\CC ^{2n}\sslash S_{n}$ becomes a morphism from
$X_{n}$ to $H_{n}$.  We have a commutative square
\begin{equation}\label{e:X-H-diagram}
\begin{CD}
X_{n}&		@>{f}>>&	\CC ^{2n}\\
@V{\rho }VV&	&		@VVV\\
H_{n}&		@>{\sigma }>>&	S^{n}\CC ^{2},
\end{CD}
\end{equation}
in which $X_{n}\subseteq H_{n}\times \CC ^{2n}$ is the set-theoretic
fiber product, with its induced reduced scheme structure.  In other
words, $X_{n}$ is the isospectral Hilbert scheme, as defined
in \cite{Hai01}.  We again write $P = \rho _{*}\Ocal _{X_{n}}$, as we
did above for a general $V\sslash G$.  Now we regard $P$ as a bundle
on $H_{n}$ rather than on $\CC ^{2n}\sslash S_{n}$.  Thus $H_{n}$ has
two different ``tautological'' bundles, the usual one $B$ and the {\it
unusual} one $P$.  The unusual tautological bundle $P$ has rank $n!$,
with an $S_{n}$ action affording the regular representation on every
fiber.  Our notation for the various schemes, bundles and morphisms
just described is identical to that in \cite{Hai01}.


The two-dimensional torus group
\begin{equation}\label{e:T2}
\TT ^{2} = (\CC ^{*})^{2}
\end{equation}
acts linearly on $\CC ^{2}$ as the group of $2\times 2$ diagonal
matrices.  We write 
\begin{equation}\label{e:tau-tq}
\tau _{t,q} = \begin{bmatrix}
t^{-1}&	0\\
0&	q^{-1}
\end{bmatrix}
\end{equation}
for its elements.  Note that when a group $G$ acts on a scheme $V$,
elements $g\in G$ act on regular functions $f\in \Ocal (V)$ as $g f =
f\circ g^{-1}$.  The inverse signs in \eqref{e:tau-tq} serve to make
$\TT ^{2}$ act on the coordinate ring $\CC [x,y]$ of $\CC ^{2}$ by the
convenient rule
\begin{equation}\label{e:T2onC[x,y]}
\tau _{t,q}x = tx;\quad \tau _{t,q}y = qy.
\end{equation}
The action of $\TT ^{2}$ on $\CC ^{2}$ induces an action on the
Hilbert scheme $H_{n}$ and all other schemes under consideration.  In
particular, $\TT ^{2}$ acts on the universal family $F_{n}$ and the
isospectral Hilbert scheme $X_{n}$, so that the projections $\pi
\colon F_{n}\rightarrow H_{n}$ and $\sigma \colon X_{n}\rightarrow
H_{n}$ are equivariant.  Hence $\TT ^{2}$ acts equivariantly on the
vector bundles $B$ and $P$.  There are induced $\TT ^{2}$ actions on
various algebraic spaces, such as the coordinate ring $\CC [\xx ,\yy
]$ of $\CC ^{2n}$, the space of global sections of any $\TT
^{2}$-equivariant vector bundle, or the fiber of such a bundle at a
torus-fixed point in $H_{n}$.  In these spaces, the $\TT ^{2}$ action
is equivalently described by a $\ZZ ^{2}$-grading.  Namely, an element
$f$ is homogeneous of degree $(r,s)$ if and only if it is a
simultaneous eigenvector of the $\TT ^{2}$ action with weight $\tau
_{t,q} f = t^{r}q^{s}f$.  Where there is an obvious natural double
grading, as in $\CC [\xx ,\yy ]$, it coincides with the weight grading
for the torus action.


We have defined the bundles whose tensor products will be the subject
of our vanishing theorems.  The theorems also describe their spaces of
global sections.  To identify these spaces, we first need to recall
the definition of the {\it polygraph} $Z(n,l)$ from \cite{Hai01}.
There, $Z(n,l)$ was defined as a certain union of linear subspaces in
$\CC ^{2n+2l}$, but it is better here to describe it first from a
Hilbert scheme point of view.  Let
\begin{equation}\label{e:Z=Xn*Fn^l}
W = X_{n}\times F_{n}^{l} \, / \, H_{n}
\end{equation}
be the fiber product over $H_{n}$ of $X_{n}$ with $l$ copies of the
universal family $F_{n}$.  The scheme $W$ is thus a closed subscheme
of $H_{n} \times \CC ^{2n+2l}$, since we have $X_{n}\subseteq
H_{n}\times \CC ^{2n}$ and $F_{n}\subseteq H_{n}\times \CC ^{2}$.  We
now define $Z(n,l)\subseteq \CC ^{2n+2l}$ to be the image of the
projection of $W$ on $\CC ^{2n+2l}$.

To see that this agrees with the original definition in \cite{Hai01},
let us identify the set $Z(n,l)$ more directly.  From
\eqref{e:X-H-diagram}, we see that a point of $X_{n}$ is an ordered
tuple 
\begin{equation}\label{e:point-of-Xn}
(I,P_{1},\ldots,P_{n})\in H_{n}\times \CC ^{2n}
\end{equation}
such that
$\sigma (I) = \llbracket P_{1},\ldots,P_{n}\rrbracket$.  In
particular, this implies $V(I) = \{P_{1},\ldots,P_{n} \}$ as a set.  A
point of $F$ is a pair $(I,Q)\in H_{n}\times \CC ^{2}$ such that $Q\in
V(I)$.  Hence a point of $W$ is a tuple
\begin{equation}\label{e:point-of-Z}
(I,P_{1},\ldots,P_{n},Q_{1},\ldots,Q_{l})\in H_{n}\times \CC ^{2n+2l}
\end{equation}
such that $\sigma (I) = \llbracket P_{1},\ldots,P_{n} \rrbracket$ and
$Q_{i}\in \{P_{1},\ldots,P_{n} \}\quad \text{for all $1\leq i\leq l$}$.
Projecting on $\CC ^{2n+2l}$, we see that 
\begin{equation}\label{e:Z(n,l)-simple}
Z(n,l) = \{(P_{1},\ldots,P_{n},Q_{1},\ldots,Q_{l})\in \CC ^{2n+2l}:
Q_{i}\in \{P_{1},\ldots,P_{n} \}\; \forall i \}.
\end{equation}
This is equivalent to the definition in \cite{Hai01}.
The scheme $W$ is flat over $H_{n}$ and reduced over the generic locus
(the open set in $H_{n}$ where the $P_{i}$ are all distinct).
Hence $W$ is reduced.  The set-theoretic description we have just
given of the projection of $W$ on $Z(n,l)$ therefore also describes a
morphism of schemes $W\rightarrow Z(n,l)$, in which we regard $Z(n,l)$
as a reduced closed subscheme of $\CC ^{2n+2l}$.


As in \cite{Hai01}, the coordinate ring of the polygraph $Z(n,l)$ will
be denoted $R(n,l)$.  Writing
\begin{equation}\label{e:xyab}
\xx ,\yy ,\aa ,\bb = x_{1}, y_{1}, \ldots , x_{n}, y_{n}, a_{1},
b_{1}, \ldots ,a_{l}, b_{l}
\end{equation}
for the coordinates on $\CC ^{2n+2l}$, we see that $R(n,l)$ is the
quotient of the polynomial ring $\CC [\xx ,\yy ,\aa ,\bb ]$ by a
suitable ideal $I(n,l)$.  Given a global regular function on $Z(n,l)$,
we may compose it with the projection $W\rightarrow Z(n,l)$ to get a
global regular function on $W$, which is the same thing as a global
section of ${P\otimes B^{\otimes l}}$ on $H_{n}$.  Hence we have a
canonical injective ring homomorphism
\begin{equation}\label{e:R(n,l)->H^0}
\psi \colon R(n,l)\hookrightarrow H^{0}(H_{n}, P\otimes B^{\otimes
l}).
\end{equation}

We can now state our first vanishing theorem, which will be proven in
Section \ref{s:van}.

\begin{thm}\label{t:van}
For all $l$ we have
\begin{gather}\label{e:H^i=0}
H^{i}(H_{n},P\otimes B^{\otimes l}) = 0\quad \text{for $i>0$, and}\\
\label{e:H^0=R(n,l)}
H^{0}(H_{n},P\otimes B^{\otimes l}) = R(n,l),
\end{gather}
where $R(n,l)$ is the coordinate ring of the polygraph
$Z(n,l)\subseteq \CC ^{2n+2l}$.
\end{thm}

The equal sign in \eqref{e:H^0=R(n,l)} is to be understood as
signifying that the homomorphism $\psi $ in \eqref{e:R(n,l)->H^0} is
an isomorphism.


Our second vanishing theorem is the analog of Theorem~\ref{t:van} for
the restriction of the tautological bundles to the {\it zero fiber}
$Z_{n} = \sigma ^{-1}(\{\underline{0} \})\subseteq H_{n}$.
In \cite{Hai98} we showed that the scheme-theoretic zero fiber is
reduced, so there is no ambiguity as to the scheme structure of
$Z_{n}$.  The ideal of the origin $\{\underline{0}\}\subseteq S^{n}\CC
^{2} = \CC ^{2n}/S_{n}$ is the homogeneous maximal ideal $\mfrak = \CC
[\xx ,\yy ]^{S_{n}}_{+}$ in the ring of invariants $\CC [\xx ,\yy
]^{S_{n}}$.  Pulled back to $H_{n}$ via $\sigma $, the elements of
$\mfrak $ represent global functions on $H_{n}$ that vanish on
$Z_{n}$.  The bundle ${P\otimes B^{\otimes l}}$ is a sheaf of $\Ocal
_{H_{n}}$-algebras, so we have a canonical inclusion
\begin{equation}\label{e:O(Hn)inH0(PBl)}
H^{0}(H_{n},\Ocal _{H_{n}})\subseteq H^{0}(H_{n}, P\otimes B^{\otimes l}).
\end{equation}
Our choice of
coordinates $\xx ,\yy ,\aa ,\bb $ on $Z(n,l)$ identifies $\CC [\xx
,\yy ]$ and $\CC [\xx ,\yy ]^{S_{n}}$ with subrings of
$R(n,l)$, in such a way that the diagram
\begin{equation}\label{e:C[x,y]inR(n,l)}
\begin{CD}
H^{0}(H_{n},\Ocal _{H_{n}})&\quad & \hookrightarrow &\quad &
H^{0}(H_{n}, P\otimes B^{\otimes l})\\
@A{\sigma ^{*}}AA&	& 		@A{\psi }AA\\
\CC [\xx ,\yy ]^{S_{n}}& \quad & \hookrightarrow &\quad & R(n,l)
\end{CD}
\end{equation}
commutes. It follows immediately that $\psi $ maps every element of
the ideal $\mfrak R(n,l)$ to a section of ${P\otimes B^{\otimes l}}$
that vanishes on $Z_{n}$.  Composing $\psi $ with restriction of
sections to the zero fiber, we get a well-defined homomorphism
\begin{equation}\label{e:psi-prime}
\psi _{1}\colon R(n,l)/\mfrak R(n,l)\rightarrow H^{0}(Z_{n}, P\otimes
B^{\otimes l}).
\end{equation}
{\it A priori}, $\psi _{1}$ need neither be injective nor surjective,
but according to our next theorem, it is both.


\begin{thm}\label{t:van-zero}
For all $l$ we have
\begin{gather}\label{e:H^i(Zn)=0}
H^{i}(Z_{n},P\otimes B^{\otimes l}) = 0\quad \text{for $i>0$, and}\\
\label{e:H^0(Zn)=R(n,l)/I} H^{0}(Z_{n},P\otimes B^{\otimes l}) =
R(n,l)/\mfrak R(n,l),
\end{gather}
where $R(n,l)$ is the polygraph coordinate ring and $\mfrak $ is the
homogeneous maximal ideal in the subring $\CC [\xx ,\yy
]^{S_{n}}\subseteq R(n,l)$.
\end{thm}

Again, the equal sign in \eqref{e:H^0(Zn)=R(n,l)/I} signifies that the
homomorphism $\psi _{1}$ in \eqref{e:psi-prime} is an isomorphism.

In a sense, Theorem \ref{t:van-zero} is a corollary to Theorem
\ref{t:van}.  Its proof uses an $\Ocal _{H_{n}}$-locally free
resolution of $\Ocal _{Z_{n}}$, which we now describe.  Afterwards, we
will prove that Theorem~\ref{t:van} implies Theorem~\ref{t:van-zero}.
The resolution we construct will be $\TT ^{2}$-equivariant.  To write
it down we need a bit more notation.  Let $\CC _{t}$ and $\CC _{q}$
denote the $1$-dimensional representations of $\TT ^{2}$ on which
$\tau _{t,q}\in \TT ^{2}$ acts by $t$ and $q$, respectively.  We write
\begin{equation}\label{e:OtOq}
\Ocal _{t} = \CC _{t}\otimes _{\CC }\Ocal _{H_{n}},\quad \Ocal _{q} =
\CC _{q}\otimes _{\CC } \Ocal _{H_{n}}
\end{equation}
for $\Ocal _{H_{n}}$ with its natural $\TT ^{2}$ action twisted by
these $1$-dimensional characters.  The $\TT ^{2}$-equivariant sheaves
$\Ocal_{t}$ and $\Ocal _{q}$ may be thought of as copies of $\Ocal
_{H_{n}}$ with respective degree shifts of $(1,0)$ and $(0,1)$.


There is a {\it trace homomorphism} of $\Ocal _{H_{n}}$-modules
\begin{equation}\label{e:tr}
\tr \colon B\rightarrow \Ocal _{H_{n}}
\end{equation}
defined as follows.  Let $\alpha \in B(U)$ be a section of $B$ on some
open set $U$.  Since $B$ is a sheaf of $\Ocal _{H_{n}}$-algebras and
also a vector bundle, there is a regular function $\tr (\alpha )\in
\Ocal _{H_{n}}(U)$ whose value at $I$ is the trace of multiplication
by $\alpha $ on the fiber $B(I)$.  The sheaf $B$ is a quotient of
${\CC [x,y]\otimes \Ocal _{H_{n}}}$, so it is generated by its global
sections $x^{r}y^{s}$ ({\it i.e.}, they span every fiber).  The trace
map is given on these sections by
\begin{equation}\label{e:tr=prs}
\tr (x^{r}y^{s})\; =\; p_{r,s}(\xx ,\yy ) \defeq \sum _{i=1}^{n}
x_{i}^{r} y_{i}^{s}.
\end{equation}
Here we regard the symmetric function $p_{r,s}\in \CC [\xx ,\yy
]^{S_{n}}$, called a {\it polarized power-sum}, as a global regular
function on $H_{n}$ pulled back from $S^{n}\CC ^{2}$ via the Chow
morphism.  To verify \eqref{e:tr=prs} we need only check it on points
$I$ in the generic locus, where the fiber $B(I) = \CC [x,y]/I$ is the
coordinate ring of a set of $n$ distinct points
$\{(x_{1},y_{1}),\ldots,(x_{n},y_{n}) \}\subseteq \CC ^{2}$.  There it
is clear that the eigenvalues of multiplication by $x^{r}y^{s}$ in
$B(I)$ are just $x_{1}^{r}y_{1}^{s},\ldots,x_{n}^{r}y_{n}^{s}$.  In
particular, $\frac{1}{n}\tr (1) = 1$, so
\begin{equation}\label{e:tr/n}
\frac{1}{n}\tr \colon B\rightarrow \Ocal _{H_{n}}
\end{equation}
is left inverse to the canonical inclusion $\Ocal _{H_{n}}
\hookrightarrow B.$ Thus we have a direct-sum decomposition of $\Ocal
_{H_{n}}$-module sheaves, or of vector bundles,
\begin{equation}\label{e:B=O+Bprime}
B = \Ocal _{H_{n}} \oplus B',\quad \text{where}\quad B' = \ker (\tr ).
\end{equation}
The projection of $B$ on its summand $B'$ is given by
$\operatorname{id}-\frac{1}{n}\tr $, so from \eqref{e:tr=prs}, we see
that $B'$ is generated by its global sections
\begin{equation}\label{e:B'global-sect}
x^{r}y^{s} - \frac{1}{n}p_{r,s}(\xx ,\yy ).
\end{equation}
Here we can omit the section corresponding to $r=s=0$, which is
identically zero.


Let $J$ be the sheaf of ideals in $B$ generated by the global sections
$x$ and $y$ and the subsheaf $B'$.  An alternative way to describe $J$
is as follows.  There are $\TT ^{2}$-equivariant sheaf homomorphisms
$\Ocal _{t}\rightarrow B$ and $\Ocal _{q}\rightarrow B$ sending the
generating section $1$ in $\Ocal _{t}$ and $\Ocal _{q}$ to $x$ and
$y$, respectively.  Combining these with the inclusion
$B'\hookrightarrow B$, we get a homomorphism of sheaves of $\Ocal
_{H_{n}}$-modules
\begin{equation}\label{e:nu}
\nu \colon B'\oplus \Ocal _{t}\oplus \Ocal _{q}\rightarrow B.
\end{equation}
Now composing ${1\otimes \nu} \colon {B\otimes (B'\oplus \Ocal
_{t}\oplus \Ocal _{q})}\rightarrow {B\otimes B}$ with the
multiplication map $\mu \colon {B\otimes B}\rightarrow B$, we get a
homomorphism of sheaves of $B$-modules
\begin{equation}\label{e:xi}
\xi \colon B\otimes (B'\oplus \Ocal _{t}\oplus \Ocal _{q})\rightarrow
B,
\end{equation}
whose image is exactly $J$.  Note that since $x$ and $y$ generate $B$
as a sheaf of $\Ocal _{H_{n}}$-algebras, the canonical homomorphism
$\Ocal _{H_{n}}\rightarrow B/J$ is surjective.  Thus $B/J$ is
identified with a quotient of $\Ocal _{H_{n}}$, which turns out to be
$\Ocal _{Z_{n}}$.


\begin{prop}\label{p:O_zero}
Let $J$ be the sheaf of ideals in $B$ generated by $x$, $y$ and $B'$.
Then $B/J$ is isomorphic as a sheaf of $\Ocal _{H_{n}}$-algebras to
$\Ocal _{Z_{n}}$.
\end{prop}

Let us recall the proof from \cite{Hai98,Hai99}, skipping some
details.  Denote by $\tilde{Z}_n$ the set-theoretic preimage $\pi
^{-1}(Z_{n})$, regarded as a reduced closed subscheme of the universal
family $F_{n}$.  Clearly the regular functions $x$, $y$ and
$p_{r,s}(\xx ,\yy )$ for $r+s>0$ vanish on $\tilde{Z}_n$.  By an old
theorem of Weyl \cite{Wey39}, the $p_{r,s}$ generate $\CC [\xx ,\yy
]^{S_{n}}$, so their vanishing defines $Z_{n}$ as a subscheme of
$H_{n}$.  Hence $\tilde{Z}_n$ is defined set-theoretically by the
vanishing of $x$, $y$ and the $p_{r,s}$, or equivalently of $x$, $y$,
and every $x^{r}y^{s}-\frac{1}{n}p_{r,s}$.  But these sections
generate $J$, so the subscheme of $F_{n}$ defined by the ideal sheaf
$J\subseteq B$ coincides set-theoretically with $\tilde{Z}_n$.  We
already know that $B/J\cong \Ocal _{Z}$ for some subscheme $Z\subseteq
H_{n}$, and this shows that $Z$ coincides set-theoretically with
$Z_{n}$.  Now $F_{n}$ is flat and finite over the non-singular scheme
$H_{n}$, hence Cohen-Macaulay.  Since $\tilde{Z}_n$ projects
bijectively on $Z_{n}$, it has codimension $n+1$ in $F_{n}$.  But
$B'\oplus \Ocal _{t}\oplus \Ocal _{q}$ is locally free of rank $n+1$,
so $J$ is everywhere locally generated by $n+1$ elements.  It follows
that $\Spec B/J$ is a local complete intersection in $F_{n}$.
Finally, one shows that $\Spec B/J$ is generically reduced, hence
reduced, which implies $B/J\cong \Ocal _{Z_{n}}$.


Our motive in reviewing this is to note that $J$ is locally a complete
intersection ideal in $B$ generated by the image under $\xi $ of any
local basis of ${B'\oplus \Ocal _{t}\oplus \Ocal _{q}}$.  Hence the
Koszul complex on the map $\xi $ in \eqref{e:xi} is a resolution of
$B/J\cong \Ocal _{Z_{n}}$.  Since everything in the construction is
$\TT ^{2}$-equivariant we deduce the following result.

\begin{prop}\label{p:res} We have a $\TT ^{2}$-equivariant locally
$\Ocal _{H_{n}}$-free resolution
\begin{equation}\label{e:res}
\cdots \rightarrow B\otimes \ext ^{k}(B'\oplus \Ocal _{t}\oplus \Ocal
_{q})\rightarrow \cdots \rightarrow B\otimes (B'\oplus \Ocal
_{t}\oplus \Ocal _{q})\underset{\xi }{\rightarrow }B\rightarrow \Ocal
_{Z_{n}}\rightarrow 0,
\end{equation}
where $\xi $ is the sheaf homomorphism in \eqref{e:xi}.
\end{prop}

As in \cite{Hai98}, it follows as a corollary that the
scheme-theoretic zero fiber is equal to the reduced zero fiber, and
that it is Cohen-Macaulay.


\begin{proof}[Proof that Theorem~\ref{t:van} implies
Theorem~\ref{t:van-zero}] Let $V\subdot $ denote the complex in
\eqref{e:res} with the final term $\Ocal _{Z_{n}}$ deleted.  The fact
that \eqref{e:res} is a resolution means that $V\subdot $ and $\Ocal
_{Z_{n}}$ are isomorphic as objects in the derived category
$D(H_{n})$.  Here and below we work in the derived category of
complexes of sheaves of $\Ocal _{H_{n}}$-modules with bounded,
coherent cohomology.  Note that $V\subdot $ is a complex of locally
free sheaves, each of which is a sum of direct summands of tensor
powers of $B$.  It follows from Theorem~\ref{t:van} that
${P\otimes B^{\otimes l}\otimes V\subdot }$ is a complex of acyclic
objects for the global section functor $\Gamma $ on $H_{n}$, so we
have
\begin{equation}\label{e:RGamma}
R\Gamma (P\otimes B^{\otimes l}\otimes V\subdot ) = \Gamma (P\otimes
B^{\otimes l}\otimes V\subdot ).
\end{equation}
%
% first two $..$ on next line have no danger of breaking badly,
% and leaving them open allows better fine adjustment
%
Now $P\otimes B^{\otimes l}\otimes V\subdot $ is isomorphic to
$P\otimes B^{\otimes l}\otimes \Ocal _{Z_{n}}$ in $D(H_{n})$, so
$H^{i}(Z_{n},{P\otimes B^{\otimes l}}) = R^{i}\Gamma ({P\otimes
B^{\otimes l}\otimes V\subdot })$ is the $i$-th cohomology of the
complex in \eqref{e:RGamma}.  This complex is zero in positive
degrees, so we deduce that $H^{i}(Z_{n},{P\otimes B^{\otimes l}}) = 0$
for $i>0$, which is the first part of Theorem~\ref{t:van-zero}.  This
is just the standard argument for the higher cohomology vanishing of a
sheaf with an acyclic left resolution.  Since $H^{i}(Z_{n},{P\otimes
B^{\otimes l}})$ is zero in negative degrees, we also deduce that the
complex in \eqref{e:RGamma} is a resolution of $H^{0}(Z_{n},{P\otimes
B^{\otimes l}})$.


Consider the last terms in this resolution:
\begin{equation}\label{e:tail}
\Gamma (P\otimes B^{\otimes l+1}\otimes (B'\oplus \Ocal _{t}\oplus \Ocal
_{q})) \underset{\Gamma (1\otimes \xi )}{\rightarrow } R(n,l+1)
\rightarrow H^{0}(Z_{n},{P\otimes B^{\otimes l}})\rightarrow 0.
\end{equation}
Here we have identified $\Gamma ({P\otimes B^{\otimes l}}\otimes B)$
with $R(n,l+1)$ using Theorem~\ref{t:van}.  To keep the notation
consistent, we denote the coordinates in $R(n,l+1)$ corresponding to
the tensor factor $B$ coming from $V\subdot $ by $x,y$ instead of the
usual $a_{l+1},b_{l+1}$.  The subring of $R(n,l+1)$ generated by the
remaining coordinates $\xx ,\yy ,\aa ,\bb $ is just $R(n,l)$, since
the projection of $Z(n,l+1)$ on these coordinates is $Z(n,l)$.  The
homomorphism $R(n,l+1)\rightarrow H^{0}(Z_{n},{P\otimes B^{\otimes
l}})$ sends $x$ and $y$ to zero and coincides on $R(n,l)$ with $\psi
_{1}$ composed with the canonical map $R(n,l)\rightarrow R(n,l)/\mfrak
R(n,l)$.  Using Theorem~\ref{t:van} we can also identify $\Gamma
({P\otimes B^{\otimes l+1}\otimes (B'\oplus \Ocal _{t}\oplus \Ocal
_{q})})$ with
\begin{equation}\label{e:R+R+R}
R(n,l+2)'\oplus R(n,l+1)\oplus R(n,l+1),
\end{equation}
where $R(n,l+2)'$ is the direct summand $\Gamma ({P\otimes B^{\otimes
l+1}\otimes B'})$ of $R(n,l+2) = \Gamma ({P\otimes B^{\otimes
l+1}\otimes B})$.  In $R(n,l+2)$ we write $x$, $y$, $x'$, $y'$ for
$a_{l+1}$, $b_{l+1}$, $a_{l+2}$, $b_{l+2}$.  By
\eqref{e:B'global-sect}, $R(n,l+2)'$ is the $R(n,l+1)$-submodule of
$R(n,l+2)$ generated by all
\begin{equation}\label{e:R'inR}
(x')^{r}(y')^{s} - \frac{1}{n}p_{r,s}(\xx ,\yy ).
\end{equation}
More precisely, $R(n,l+2)$ is generated as an $R(n,l+1)$-module by the
monomials $(x')^{r}(y')^{s}$, and the projection on the summand
$R(n,l+2)'$ is the homomorphism of $R(n,l+1)$ modules mapping
$(x')^{r}(y')^{s}$ to the expression in \eqref{e:R'inR}.  Although we
are implicitly relying on Theorem~\ref{t:van} to guarantee that this
is well defined, it can also be shown directly.


The map $\Gamma ({1\otimes \xi })$ in \eqref{e:tail} now becomes the
$R(n,l+1)$-module homomorphism
\begin{equation}\label{e:Gam(1.xi)}
R(n,l+2)'\oplus R(n,l+1)\oplus R(n,l+1)\rightarrow R(n,l+1)
\end{equation}
given on the first summand by $(x')^{r}(y')^{s}\mapsto x^{r}y^{s}$ and
on the second and third summands by multiplication by $x$ and $y$,
respectively.  Its image is therefore the ideal in $R(n,l+1)$
generated by $x$, $y$ and all $x^{r}y^{s}-\frac{1}{n}p_{r,s}(\xx ,\yy
)$, or equivalently, the ideal
\begin{equation}\label{e:J}
J = (x,y)+\mfrak R(n,l+1).
\end{equation} 
Since $x$ and $y$ generate $R(n,l+1)$ as an $R(n,l)$-module, the
inclusion $R(n,l)\subseteq R(n,l+1)$ induces a surjective ring
homomorphism
\begin{equation}\label{e:R->R/J}
R(n,l)\rightarrow R(n,l+1)/J
\end{equation}
with kernel 
\begin{equation}\label{e:I}
I = R(n,l)\cap J.
\end{equation}
By \eqref{e:tail}, we have $R(n,l)/I\cong R(n,l+1)/J\cong
H^{0}(Z_{n},{P\otimes B^{\otimes l}})$.  The isomorphism here is
induced by $\psi _{1}$.  Thus it only remains to show that $I = \mfrak
R(n,l)$.


Clearly, $I$ contains $\mfrak R(n,l)$, so we are to show
that the homomorphism
\begin{equation}\label{e:R/mR->R/I}
\zeta \colon R(n,l)/\mfrak R(n,l)\rightarrow R(n,l)/I \cong R(n,l+1)/J
\end{equation}
is injective.  For this we construct its left inverse.  From the
equation $R(n,l+1) = \Gamma ({P\otimes B^{\otimes l+1}})$ and the
decomposition $B = \Ocal _{H_{n}}\oplus B'$, taken in the last tensor
factor of $B^{\otimes l+1}$, we see that $R(n,l)$ is a direct summand
of $R(n,l+1)$ as an $R(n,l)$-module.  Using \eqref{e:tr=prs} and
\eqref{e:tr/n}, we obtain the formula
\begin{equation}\label{e:theta}
\theta (x^{r}y^{s}) = \frac{1}{n}p_{r,s}(\xx ,\yy )
\end{equation}
for the projection $\theta \colon R(n,l+1)\rightarrow R(n,l)$.  Now,
$\theta $ is a homomorphism of $R(n,l)$-modules and $\mfrak $ is a
subset of $R(n,l)$, so $\theta $ carries $\mfrak R(n,l+1)$ into
$\mfrak R(n,l)$.  The monomials $x^{r}y^{s}$ with $r+s>0$ generate
$(x,y)R(n,l+1)$ as an $R(n,l)$-module, so \eqref{e:theta} shows that
$\theta $ also carries $(x,y)R(n,l+1)$ into $\mfrak R(n,l)$.  Hence
$\theta $ induces a map
\begin{equation}\label{e:thetabar}
\overline{\theta }\colon R(n,l+1)/J\rightarrow R(n,l)/\mfrak R(n,l).
\end{equation}
The endomorphism $\overline{\theta }\circ \zeta $ of $R(n,l)/\mfrak
R(n,l)$ is a homomorphism of $R(n,l)$-modules, so it is the identity,
and $\overline{\theta }$ is the required left inverse of $\zeta $.
\end{proof}


\section{Character formulas} \label{s:char}

Theorems~\ref{t:van} and \ref{t:van-zero} allow us to identify the
ring of diagonal coinvariants and the polygraph coordinate ring
$R(n,l)$, among other things, with spaces of global sections of $\TT
^{2}$-equivariant coherent sheaves on $H_{n}$.  When the higher
cohomology vanishes, we can calculate the $\TT ^{2}$ character of the
space of global sections, or what is the same, its Hilbert series as a
doubly graded module, using the Lefschetz formula of Atiyah and
Bott \cite{AtBo66}.  We will apply this method to obtain explicit
character formulas for spaces of interest in the Hilbert scheme
context.  As we shall see, the resulting formulas are naturally
expressed in terms of operators arising in the theory of Macdonald
polynomials.

Let $M = \bigoplus  M_{r,s}$ be a finitely-generated doubly graded
module over $\CC [\xx ,\yy ]$ or $\CC [\xx ,\yy ]^{S_{n}}$.  The {\it
Hilbert series} of $M$ is the Laurent series in two variables
\begin{equation}\label{e:HM}
\Hcal _{M}(q,t) = \sum _{r,s} t^{r}q^{s}\dim (M_{r,s}).
\end{equation}
If $M$ is finite-dimensional as a vector space over $\CC $, then
$\Hcal _{M}(q,t) = \tr (M,\tau _{t,q})$ is the character of $M$ as a
$\TT ^{2}$-module in the strict sense.  In general, it is a good idea
to think of $\Hcal _{M}(q,t)$ as a formal $\TT ^{2}$ character, for
reasons that will become apparent below.  The Laurent series $\Hcal
_{M}(q,t)$ is a rational function of $q$ and $t$.  When $M$ is a $\CC
[\xx ,\yy ]$-module this is well-known and can be shown easily by
calculating the Hilbert series using a finite graded free resolution
of $M$.  When $M$ is a $\CC [\xx ,\yy ]^{S_{n}}$-module, one obtains
the same result by regarding $M$ as a module over $\CC [p_{1}(\xx
),p_{1}(\yy ), \ldots, p_{n}(\xx ),p_{n}(\yy )]$, since the power sums
$p_{k}(\xx )$, $p_{k}(\yy )$ form a doubly homogeneous system of
parameters in $\CC [\xx ,\yy ]^{S_{n}}$.  Now let $A$ be a $\TT
^{2}$-equivariant coherent sheaf on $H_{n}$.  The Chow morphism
$\sigma \colon H_{n}\rightarrow S^{n}\CC ^{2}$ is projective, and
$S^{n}\CC ^{2}$ is affine, so the sheaf cohomology modules
$H^{i}(H_{n},A)$ are finitely-generated $\TT ^{2}$-equivariant---which
is to say, doubly graded---$\CC [\xx ,\yy ]^{S_{n}}$-modules.  We
denote their Hilbert series by
\begin{equation}\label{e:HiA}
\Hcal ^{i}_{A}(q,t) = \Hcal _{H^{i}(H_{n},A)}(q,t).
\end{equation}


The Atiyah--Bott formula expresses the Euler characteristic
\begin{equation}\label{e:chiA}
\chi _{A}(q,t) \defeq \sum _{i}(-1)^{i}\Hcal ^{i}_{A}(q,t)
\end{equation}
as a sum of local contributions from the $\TT ^{2}$-fixed points of
$H_{n}$.  These local contributions are described by data associated
with partitions of $n$.  Let us fix some notation.  We write a
partition of $n$ as $\mu =(\mu _{1}\geq \mu _{2}\geq \cdots \geq \mu
_{l}>0)$, with the understanding that $\mu _{i}=0$ for $i>l$.  The
{\it Ferrers diagram} of $\mu $ is the set of lattice points
\begin{equation}\label{e:d(mu)}
d(\mu ) = \{(i,j)\in \NN \times \NN : j <\mu _{i+1} \}.
\end{equation}
By convention, the diagram is drawn with the $i$-axis vertical and the
$j$-axis horizontal, so the parts of $\mu $ are the lengths of the
rows, and $(0,0)$ is the lower left corner.  The {\it arm} $a(x)$ and
{\it leg} $l(x)$ of a point $x\in d(\mu )$ denote the number of points
strictly to the right of $x$ and above $x$, respectively, as indicated
in this example:
\begin{equation}\label{e:arm-leg-pix}
\mu =(5,5,4,3,1)\qquad 
\begin{array}[c]{cccccc}
\bullet &	\hbox to 0pt{\hss $\scriptstyle l(x)$\hss }\\
\cline{2-2}
\bullet &	\multicolumn{1}{|c|}{\bullet }&	\bullet \\
\bullet &	\multicolumn{1}{|c|}{\bullet }&	\bullet &	\bullet \\
\cline{2-5}
\bullet &	\multicolumn{1}{|c|}{\llap{${}_{x}$}\bullet } &	\bullet &
    \bullet&	 \multicolumn{1}{c|}{\bullet }&	{\scriptstyle a(x)} \\
\cline{2-5}
\llap{${}_{(0,0)}$} \bullet &	\bullet &	\bullet &
\bullet &	\bullet  
\end{array}
\qquad a(x) = 3,\quad l(x) = 2.
\end{equation}
To each partition $\mu $ is associated a monomial ideal
\begin{equation}\label{e:Imu}
I_{\mu } = \CC \cdot \{x^{r}y^{s}:(r,s)\not \in d(\mu ) \}\subseteq
\CC [x,y].
\end{equation}
A $\CC $-basis of 
 $\CC [x,y]/I_{\mu }$ is given by the set of monomials not in $I_{\mu }$,
\begin{equation}\label{e:Bmu-set}
\Bcal _{\mu } =  \{x^{r}y^{s}:(r,s) \in d(\mu ) \}.
\end{equation}
In particular, $\dim _{\CC }\CC [x,y]/I_{\mu } = n$, so $I_{\mu }$ is
a point of $H_{n}$.


\begin{prop}\label{p:fixed-pts}
The $\TT ^{2}$-fixed points of $H_{n}$ are the ideals $I_{\mu }$ for
all partitions $\mu $ of $n$.  The cotangent space of $H_{n}$ at
$I_{\mu }$ has a basis of $\TT ^{2}$-eigenvectors $\{u_{x},d_{x}:x\in
d(\mu ) \}$ with eigenvalues
\begin{equation}\label{e:eigens-at-mu}
\tau _{t,q}d_{x} = t^{1+l(x)}q^{-a(x)}d_{x},\quad \tau _{t,q}u_{x} =
t^{-l(x)}q^{1+a(x)}u_{x}.
\end{equation}
\end{prop}

\begin{proof}
An ideal $I\subseteq \CC [x,y]$ is $\TT ^{2}$-fixed if and only if it
is doubly homogeneous, or equivalently, a monomial ideal.  This
establishes the first part.  The eigenvalues (expressed somewhat
differently) were determined by Ellingsrud and
Str\"omme \cite{ElSt87}.  The basis elements $u_{x}$, $d_{x}$ are
given explicitly in terms of local coordinates
in \cite[Corollary~2.5]{Hai98}.
\end{proof}

Now we give the Atiyah--Bott formula as it applies in our context.
For simplicity we only state it for vector bundles, {\it i.e.},
locally free sheaves, which is all we need.

\begin{prop}\label{p:AB}
Let $A$ be a $\TT ^{2}$-equivariant locally free sheaf of finite rank
on $H_{n}$.  Then
\begin{equation}\label{e:AB}
\chi _{A}(q,t) = \sum _{|\mu | = n} \frac{ \Hcal _{A(I_{\mu})}(q,t) }{
\prod _{x\in d(\mu )} (1-t^{1+l(x)}q^{-a(x)}) (1-t^{-l(x)}q^{1+a(x)})
}.
\end{equation}
\end{prop}

\begin{proof}
What we have written is the classical formula in Theorem~2 of
\cite{AtBo66}, evaluated on the data in Proposition~\ref{p:fixed-pts}.
Since $H_{n}$ is not a projective variety, however, and the left-hand
side in \eqref{e:AB} is only a formal $\TT ^{2}$ character, some
further justification is required.  Various authors have extended the
classical formula to more general contexts and given algebraic proofs.
We will use the following corollary to a very general theorem of
Thomason \cite[Th\'eor\`eme~3.5]{Thom92}.


\begin{prop}\label{p:Thomason}
Let $T = \TT ^{d} = \Spec \CC
[t_{1},t_{1}^{-1},\ldots,t_{d},t_{d}^{-1}]$ be an algebraic torus, $X$
and $Y$ separated schemes of finite type over $\CC $ on which $T$
acts, and $f\colon X\rightarrow Y$ a $T$-equivariant proper morphism.
Assume $X$ is non-singular.  Let $K_{0}(T,X)$, $K_{0}(T,Y)$, {\it
etc.}\ denote the Grothendieck groups of $T$-equivariant coherent
sheaves, and $K^{0}(T,X)$, {\it etc.}\ the Grothendieck rings of
$T$-equivariant algebraic vector bundles.  Recall that (for any $X$)
$K_{0}(T,X)$ is a $K^{0}(T,X)$-module and $K^{0}(T,X)$ is an algebra
over the representation ring $R(T)$, which we identify with $\ZZ
[\boldt ,\boldt ^{-1}] = \ZZ
[t_{1},t_{1}^{-1},\ldots,t_{d},t_{d}^{-1}]$.  Define
\begin{equation}\label{e:K0(T,X)(0)}
K_{0}(T,X)_{(0)} = \QQ (\boldt )\otimes _{\ZZ [\boldt ,\boldt ^{-1}]}
K_{0}(T,X),
\end{equation}
and similarly
for $Y$, {\it etc.}.  Then the following hold.

(1) Let $N$ be the conormal bundle of the fixed-point locus $X^{T}$ in
$X$, and set $\ext N = \sum _{i}(-1)^{i}[\ext ^{i}N]\in
K^{0}(T,X^{T})$.  Then $\ext N$ is invertible in
$K^{0}(T,X^{T})_{(0)}$.

(2) Let $f_{*}\colon K_{0}(T,X)_{(0)}\rightarrow K_{0}(T,Y)_{(0)}$
be the homomorphism induced by the derived pushforward, that is,
$f_{*}[A] = \sum _{i}(-1)^{i}[R^{i}f_{*}A]$, and let $f_{*}^{T}\colon
K_{0}(T,X^{T})_{(0)}\rightarrow K_{0}(T,Y^{T})_{(0)}$ denote the same
for the fixed-point loci.  Then
\begin{equation}\label{e:Thomason}
f_{*}[A] = i_{*}f_{*}^{T}\left( (\ext N)^{-1}\cdot \sum
_{k}(-1)^{k}[\operatorname{Tor}_{k}^{\Ocal _{X}}(\Ocal _{X^{T}}, A)]
\right), 
\end{equation}
where $i_{*}\colon K_{0}(T,Y^{T})_{(0)}\rightarrow K_{0}(T,Y)_{(0)}$
is induced by $i\colon Y^{T}\hookrightarrow Y$.
\end{prop}


To obtain \eqref{e:AB}, we apply Thomason's theorem with $T = \TT
^{2}$ and $f\colon X\rightarrow Y$ the Chow morphism $\sigma \colon
H_{n}\rightarrow S^{n}\CC ^{2}$.  The group $K_{0}(T,Y)$ is identified
with the Grothendieck group of finitely-generated doubly graded $\CC
[\xx ,\yy ]^{S_{n}}$-modules.  The Hilbert series $\Hcal _{M}(q,t)$
only depends on the class $[M]\in K_{0}(T,Y)$ of $M$, and so induces a
$\ZZ [q,q^{-1},t,t^{-1}]$-linear map
\begin{equation}\label{e:H-on-K0}
\Hcal \colon K_{0}(T,Y)_{(0)}\rightarrow \QQ (q,t).
\end{equation}
The fixed-point locus $Y^{T}$ is a point, so $K_{0}(T,Y^{T})_{(0)} =
\QQ (q,t)$, and $\Hcal \circ i_{*}$ is the identity map on $\QQ
(q,t)$.  Similarly, $X^{T}$ is the finite set $\{I_{\mu }:|\mu |=n \}$
and $K_{0}(T,X^{T})_{(0)}$ is the direct sum of copies of $\QQ (q,t)$,
one for each $\mu $.  With these identifications, $f_{*}^{T}$ is just
summation over $\mu $.  Applying $\Hcal $ to both sides in
\eqref{e:Thomason} yields \eqref{e:AB}.
\end{proof}


Some of our sheaves and spaces have $S_{n}$ actions, so we need to
sharpen our notation a bit to keep track of them.  Recall that the
{\it Frobenius characteristic map} from $S_{n}$ characters to
symmetric functions is defined by
\begin{equation}\label{e:Fmap}
\phi \chi = \frac{1}{n!}\sum _{w\in S_{n}}\chi (w) p_{\tau (w)}(z),
\end{equation}
where $\tau (w)$ is the partition of $n$ given by the cycle lengths in
the expression for the permutation $w$ as a product of disjoint
cycles, and $p_{\lambda }(z) = p_{\lambda _{1}}(z)\cdots p_{\lambda
_{l}}(z)$ denotes the power-sum symmetric function.  The irreducible
characters of $S_{n}$ are then given by the identity
\begin{equation}\label{e:Fchi-lambda}
\phi \chi ^{\lambda } = s_{\lambda }(z),
\end{equation}
where $s_{\lambda }(z)$ is a Schur function.  Here and below we always
work in the algebra
\begin{equation}\label{e:sym-funcs}
\Lambda = \Lambda _{\QQ (q,t)}(z)
\end{equation}
of symmetric functions in infinitely many variables $z =
z_{1},z_{2},\ldots$ with coefficients in $\QQ (q,t)$.  As $\lambda $
runs over partitions of $n$, the power-sums $p_{\lambda }(z)$, Schur
functions $s_{\lambda }(z)$, Macdonald polynomials $P_{\lambda
}(z;q,t)$, and so forth are bases of the homogeneous subspace $\Lambda
_{n}$ of degree $n$ in $\Lambda $.  Occasionally below we will use
{\it plethystic substitution}, also known as $\lambda $-ring notation.
Let $A$ be an algebra of polynomials or formal series in some alphabet
of indeterminates, with coefficients in $\QQ (q,t)$.  Given $Y\in A$,
we define $p_{k}[Y]$ to be the result of replacing each indeterminate
in $Y$, including $q$ and $t$, with its $k$-th power.  The algebra
$\Lambda $ is freely generated over $\QQ (q,t)$ by the power-sums
$p_{k}(z)$, so there is a unique $\QQ (q,t)$-linear homomorphism
\begin{equation}\label{e:evY}
\ev _{Y}\colon \Lambda \rightarrow A,\quad \ev _{Y}p_{k}(z) = p_{k}[Y].
\end{equation}
We now define for all $f\in \Lambda $, $Y\in A$:
\begin{equation}\label{e:f[Y]}
f[Y] = \ev _{Y}f(z).
\end{equation}
We will specifically need the following instances of this construction.
\begin{itemize}
\item Setting (here and throughout) $Z = z_{1}+z_{2}+\cdots $, we recover
$f(z) = f[Z]$.
\item $f \left[\frac{Z}{1-t} \right]$ is the image of $f$ under the
automorphism of $\Lambda $ sending $p_{k}(z)$ to $p_{k}(z)/(1-t^{k})$.
We can equate $f \left[\frac{Z}{1-t} \right]$ with
$f(z,tz,t^{2}z,\ldots)$, provided we interpret the coefficients of the
latter expression, which are rational Laurent series in $t$, as
rational functions.  The same holds with $q$ in place of $t$.
\item If $Y = a_{1}+\cdots +a_{k}$ is a sum of monomials $a_{i}$ in
the indeterminates, each with coefficient $1$, then $f[Y] =
f(a_{1},\ldots,a_{k})$.
\end{itemize}


Now let $M$ be a finitely-generated doubly graded $\CC [\xx ,\yy
]^{S_{n}}$-module with an $S_{n}$ action that respects the grading,
{\it i.e.}, it commutes with the $\TT ^{2}$ action.  For instance, $M$
might be a doubly graded $\CC [\xx ,\yy ]$ module with an equivariant
$S_{n}$ action, regarded as a $\CC [\xx ,\yy ]^{S_{n}}$-module.  We
denote by $V^{\lambda }$ the irreducible representation of $S_{n}$
with character $\chi ^{\lambda }$.  Then $M$ has a canonical
direct-sum decomposition
\begin{equation}\label{e:isotypic}
M = \bigoplus _{|\lambda | = n}V^{\lambda }\otimes M_{\lambda },\quad
M_{\lambda } \defeq \Hom ^{S_{n}}(V^{\lambda },M),
\end{equation}
in which each $M_{\lambda }$ is a doubly graded $\CC [\xx ,\yy
]^{S_{n}}$-module.   We define the {\it Frobenius series} of $M$ to be
\begin{equation}\label{e:Frobseries}
\Fcal _{M}(z;q,t) \defeq  \sum _{|\lambda | = n}\Hcal _{M_{\lambda
}}(q,t)s_{\lambda }(z) = \sum _{r,s}t^{r}q^{s}\phi \ch (M_{r,s}).
\end{equation}
The second equality follows from \eqref{e:Fchi-lambda} and shows that
the Frobenius series is a generating function for the characters $\ch
(M_{r,s})$ in the same way that the Hilbert series is a generating
function for the dimensions.  The Hilbert series can be recovered from
the Frobenius series by the formula
\begin{equation}\label{e:HfromF}
\Hcal _{M}(q,t) = \langle s_{1}^{n},\Fcal _{M}(z;q,t) \rangle,
\end{equation}
where $\langle \cdot ,\cdot  \rangle$ is the usual Hall inner
product on symmetric functions.


If $A$ is a $\TT ^{2}$-equivariant coherent sheaf on $H_{n}$ with an
$S_{n}$ action commuting with the $\TT ^{2}$ action, then $A$ has a
decomposition 
\begin{equation}\label{e:Aisotypic}
A = \bigoplus _{|\lambda |=n}{V^{\lambda }\otimes_{\CC } A_{\lambda }}
\end{equation}
as in \eqref{e:isotypic}, inducing the decomposition
\eqref{e:isotypic} for the cohomology modules $M = H^{i}(H_{n},A)$.
We set
\begin{equation}\label{e:FiA}
\Fcal ^{i}_{A}(z;q,t) = \Fcal _{H^{i}(H_{n},A)}(z;q,t),\quad \chi
\Fcal _{A}(z;q,t) = \sum _{i}(-1)^{i}\Fcal ^{i}_{A}(q,t).
\end{equation}
Then for $A$ locally free we immediately obtain the Frobenius series
version of \eqref{e:AB}:
\begin{equation}\label{e:FAB}
\chi \Fcal _{A}(q,t) = \sum _{|\mu | = n} \frac{ \Fcal
_{A(I_{\mu})}(q,t) }{ \prod _{x\in d(\mu )} (1 - t^{1+l(x)}q^{-a(x)})
(1 - t^{-l(x)}q^{1+a(x)}) }.
\end{equation}
Let us now evaluate this in some specific cases.


\subsection*{Character formula for $R(n,l)$}

Taking $A = {P\otimes B^{\otimes l}}$, the $S_{n}$ action on $A$
is induced by that on $P$.  By Theorem~\ref{t:van}, we have
\begin{equation}\label{e:chi=F-R(n,l)}
\Fcal _{R(n,l)}(z;q,t) = \chi \Fcal _{A}(z;q,t).
\end{equation}
To calculate this using \eqref{e:FAB} we must evaluate
\begin{equation}\label{e:PBl-fiber}
\Fcal _{(P\otimes B^{\otimes l})(I_{\mu })}(z;q,t) = \Fcal _{P(I_{\mu
})}(z;q,t)\Hcal _{B(I_{\mu })}(q,t)^{l}.
\end{equation}
The set $\Bcal _{\mu }$ in \eqref{e:Bmu-set} is a doubly homogeneous
basis of $B(I_{\mu }) = \CC [x,y]/I_{\mu }$, so we have
\begin{equation}\label{e:Bmu(q,t)}
\Hcal _{B(I_{\mu })}(q,t) = B_{\mu }(q,t) \defeq \sum _{(r,s)\in d(\mu
)}t^{r}q^{s}.
\end{equation}
The Frobenius series of $P(I_{\mu })$ is given by the {\it  transformed
Macdonald polynomial}
\begin{equation}\label{e:Htild}
\Htild _{\mu }(z;q,t) \defeq t^{n(\mu )} J_{\mu }\left[
\smallfrac{Z}{1-t^{-1}}; q,t^{-1} \right],
\end{equation}
where $J_{\mu }$ is the {\it integral form} Macdonald polynomial
defined in \cite[VI, eq.~(8.3)]{Mac95}, and $n(\mu ) = \sum _{i}(i-1)\mu
_{i}$.  Equivalently,
\begin{equation}\label{}
\Htild _{\mu }(z;q,t) = \sum _{\lambda }\Ktild _{\lambda \mu }(q,t)
s_{\lambda }(z),\quad \Ktild _{\lambda \mu }(q,t) = t^{n(\mu
)}K_{\lambda \mu }(q,t^{-1}),
\end{equation}
where $K_{\lambda \mu }(q,t)$ is the Kostka--Macdonald
coefficient \cite[VI, eq.~(8.11)]{Mac95}.


\begin{prop}[\cite{Hai01}]\label{p:FP(Imu)=Hmu}
We have $\Fcal _{P(I_{\mu })}(z;q,t) = \Htild _{\mu }(z;q,t)$.
\end{prop}

We remark in passing that in the decomposition \eqref{e:Aisotypic} for
$P$, say $P = \bigoplus _{\lambda }{V^{\lambda }\otimes P_{\lambda
}}$, the {\it character bundles} $P_{\lambda }$ are analogous to the
character bundles on $\CC ^{2}\sslash G$ introduced by
Gonzalez-Sprinberg and Verdier \cite{GoSpVer83} in connection with the
classical ($G\subseteq \SL _{2}$) McKay correspondence.
Proposition~\ref{p:FP(Imu)=Hmu} identifies the Kostka-Macdonald
coefficient $\Ktild _{\lambda \mu }(q,t)$ as the Hilbert series of the
fiber $P_{\lambda }(I_{\mu })$.  In particular, this shows $\Ktild
_{\lambda \mu }(q,t)\in \NN [q,t]$, which was the main combinatorial
theorem in \cite{Hai01}.

Combining the proposition with the equations preceding it, we arrive
at the following result.

\begin{thm}\label{t:FR(n,l)}
The Frobenius series of $R(n,l)$ is given by
\begin{equation}\label{e:FR(n,l)}
\Fcal _{R(n,l)}(z;q,t) = \sum _{|\mu | = n} \frac{ B_{\mu }(q,t)^{l}
\Htild _{\mu }(z;q,t) }{ \prod _{x\in d(\mu )} (1 -
t^{1+l(x)}q^{-a(x)}) (1 - t^{-l(x)}q^{1+a(x)}) }.
\end{equation}
\end{thm}

We can express this more succinctly in terms of the linear operator
$\Delta $ on $\Lambda $ defined by
\begin{equation}\label{e:Delta}
\Delta \Htild _{\mu }(z;q,t) = B_{\mu }(q,t)\Htild _{\mu }(z;q,t).
\end{equation}
This operator was introduced in \cite{GaHa96}, where we gave a direct
plethystic expression for it \cite[Theorem~2.2]{GaHa96}.


\begin{lem}\label{l:Z/1-t}
Let $M$ be a finitely-generated doubly graded $\CC [\xx ,\yy ]$-module
with an equivariant $S_{n}$ action.  If the $\xx$ variables
$x_{1},\ldots,x_{n}$ form an $M$-regular sequence, then
\begin{equation}\label{e:Z[1-t]}
\Fcal _{M}(z;q,t) = \Fcal _{M/(\xx )M}\left[\smallfrac{Z}{1-t};q,t \right],
\end{equation}
and similarly with $\yy $ and $q$ in place of $\xx $ and $t$.
\end{lem}

\begin{proof}
For a module over a local ring, this was proven in \cite[Proposition
5.3]{Hai99}.  The same proof applies in the graded setting essentially
without change.
\end{proof}

\begin{lem}\label{l:R(n,0)}
The Frobenius series of $\CC [\xx ,\yy ]$ is given by
\begin{equation}\label{e:R(n,0)}
\Fcal _{\CC [\xx ,\yy ]}(z;q,t) = h_{n}\left[\smallfrac{Z}{(1-q)(1-t)} \right],
\end{equation}
where $h_{n}(z)$ is the complete homogeneous symmetric function of
degree $n$.
\end{lem}

\begin{proof}
Apply Lemma~\ref{l:Z/1-t} first to the regular sequence $\xx $ in $
\CC [\xx ,\yy ]$, then to $\yy $ in $\CC [\yy ]$.  This reduces
\eqref{e:R(n,0)} to $\Fcal _{\CC }(z;q,t)  = h_{n}(z) = s_{(n)}(z)$,
which is correct since $\CC $ is the trivial representation in degree
$(0,0)$.
\end{proof}


\begin{cor}\label{c:FR(n,l)}
The formula \eqref{e:FR(n,l)} in Theorem~\ref{t:FR(n,l)} is equivalent
to
\begin{equation}\label{e:FR(n,l)-delta}
\Fcal _{R(n,l)}(z;q,t) = \Delta ^{l}h_{n}\left[\smallfrac{Z}{(1-q)(1-t)} \right].
\end{equation}
\end{cor}

\begin{proof}
From \eqref{e:FR(n,l)} it is clear that $\Fcal _{R(n,l)}(z;q,t) =
\Delta ^{l}\Fcal _{R(n,0)}(z;q,t) $.  But $R(n,0) = \CC [\xx ,\yy ]$.
\end{proof}

Note that the case $l=0$ gives a geometric interpretation and proof of
one of the basic identities in the theory of Macdonald
polynomials \cite[Theorem~2.8]{GaHa96}:
\begin{equation}\label{e:hn[Z/M]}
h_{n}\left[\smallfrac{Z}{(1-q)(1-t)} \right] = \sum _{|\mu | = n}
\frac{ \Htild _{\mu }(z;q,t) }{ \prod _{x\in d(\mu )} (1 -
t^{1+l(x)}q^{-a(x)}) (1 - t^{-l(x)}q^{1+a(x)}) }.
\end{equation}
From the preceding corollary we obtain a formula for the Hilbert
series of $R(n,l)$.


\begin{cor}\label{c:HR(n,l)} We have
\begin{align}\label{e:HR(n,l)a}
\Hcal _{R(n,l)} &	 = \langle s_{1}^{n}(z),\Delta ^{l}
h_{n}\left[\smallfrac{Z}{(1-q)(1-t)} \right]\rangle\\
\label{e:HR(n,l)b}
&	 = \frac{1}{(1-q)^{n}(1-t)^{n}}\langle e_{n}(z),
\Delta ^{l}s_{1}^{n}(z) \rangle,
\end{align}
where $e_{n}(z)$ is the $n$-th elementary symmetric function.
\end{cor}

\begin{proof}
The first equation is immediate from Corollary~\ref{c:FR(n,l)}.  For
the second, recall from \cite{GaHa96} that the transformed Macdonald
polynomials are orthogonal with respect to the inner product
\begin{equation}\label{e:<>*}
\langle f,g \rangle_{*} \defeq \langle \omega f[Z(1-q)(1-t)], g \rangle,
\end{equation}
where $\omega $ is the familiar involution on $\Lambda $ defined by
$\omega e_{k}(z) = h_{k}(z)$.  Any operator with the $\Htild _{\mu
}(z;q,t)$ as eigenfunctions, and $\Delta $ in particular, is therefore
self-adjoint with respect to $\langle \cdot ,\cdot \rangle_{*}$.
Hence
\begin{equation}
\begin{aligned}
\langle s_{1}^{n}(z), \Delta ^{l}h_{n}\left[\smallfrac{Z}{(1-q)(1-t)}
\right]\rangle &   =	\langle \omega s_{1}^{n}\left[\smallfrac{Z}{(1-q)(1-t)}
\right], \Delta ^{l}h_{n}\left[\smallfrac{Z}{(1-q)(1-t)}
\right] \rangle_{*}\\
&	= \frac{1}{(1-q)^{n}(1-t)^{n}} \langle \Delta
^{l}s_{1}^{n}(z),h_{n}\left[\smallfrac{Z}{(1-q)(1-t)} \right] \rangle_{*} \\  
& = \frac{1}{(1-q)^{n}(1-t)^{n}}\langle e_{n}(z), \Delta
^{l}s_{1}^{n}(z) \rangle.
\end{aligned}
\end{equation}
\end{proof}


\subsection*{Character formula for diagonal coinvariants}

The {\it ring of coinvariants} for the diagonal action of $S_{n}$ on
$\CC ^{2n}$ is, by definition,
\begin{equation}\label{e:Rn}
R_{n} = \CC [\xx ,\yy ]/\mfrak \CC [\xx ,\yy ],
\end{equation}
where $\mfrak $ is the homogeneous maximal ideal in $\CC [\xx ,\yy
]^{S_{n}}$.  Ignoring its ring structure, $R_{n}$ is isomorphic as a
doubly graded $S_{n}$-module to the space of {\it diagonal harmonics}
\begin{equation}\label{e:DHn}
\Dh _{n} = \{f\in \CC [\xx ,\yy ]: p(\partial \xx ,\partial \yy )f =0\;
\forall p\in \mfrak \}.
\end{equation}
Its Frobenius series was the subject of a series of combinatorial
conjectures by the author and others in \cite{Hai94}.  Later,
in \cite{GaHa96}, Garsia and the author showed that these conjectures
would follow from a conjectured master formula giving $\Fcal
_{R_{n}}(z;q,t)$ in terms of Macdonald polynomials, which we will now
prove.


From Theorem~\ref{t:van-zero}, with $l=0$, we obtain
\begin{equation}\label{e:FRn=chi}
\Fcal _{R_{n}}(z;q,t)  = \chi \Fcal _{P\otimes \Ocal _{Z_{n}}}(z;q,t).
\end{equation}
To calculate this using \eqref{e:FAB}, we replace $\Ocal _{Z_{n}}$
with the resolution $V\subdot $ given by the complex in \eqref{e:res}
with the final term deleted.  This gives
\begin{equation}\label{e:chiRn}
\chi \Fcal _{P\otimes \Ocal _{Z_{n}}}(z;q,t) = \sum _{k=0}^{n+1}
(-1)^{k}\chi \Fcal _{P\otimes V_{k}}(z;q,t),
\end{equation}
where $V_{k} = B\otimes \ext ^{k}(B'\oplus \Ocal _{t}\oplus \Ocal
_{q})$.  The eigenvalues of $\tau _{t,q}\in \TT ^{2}$ on the fiber
$(B'\oplus \Ocal _{t}\oplus \Ocal _{q})(I_{\mu })$ are $q$ and $t$,
from the summand $\Ocal _{t}\oplus \Ocal _{q}$, and
$\{t^{r}q^{s}:(r,s)\in d(\mu )\setminus \{(0,0) \} \}$, from the basis
$\Bcal _{\mu }\setminus \{1 \}$ of $B'(I_{\mu })$.  The Hilbert series
of $\ext ^{k}(B'\oplus \Ocal _{t}\oplus \Ocal _{q})(I_{\mu })$ is
the $k$-th elementary symmetric function of these eigenvalues, and its
alternating sum over $k$ is therefore $(1-q)(1-t)\Pi _{\mu }(q,t)$,
where
\begin{equation}\label{e:Pi-mu}
\Pi _{\mu }(q,t) \defeq  \prod _{\substack{(r,s)\in d(\mu )\\
(r,s)\not =(0,0)}} (1-t^{r}q^{s}).
\end{equation}
Hence we have
\begin{equation}\label{e:FVk}
 \sum _{k=0}^{n+1} (-1)^{k}\Fcal _{(P\otimes V_{k})(I_{\mu })}(z;q,t)
= (1-q)(1-t)\Pi _{\mu }(q,t)B_{\mu }(q,t)\Htild _{\mu }(z;q,t),
\end{equation}
and \eqref{e:FAB} yields the following character formula for the
diagonal coinvariants.


\begin{thm}\label{t:FRn}
The Frobenius series of the coinvariant ring $R_{n}$, or of the
diagonal harmonics $\Dh _{n}$, is given by
\begin{equation}\label{e:FRn}
\Fcal _{R_{n}}(z;q,t) = \sum _{|\mu | = n} \frac{ (1-q)(1-t) \Pi _{\mu
}(q,t) B_{\mu }(q,t) \Htild _{\mu }(z;q,t) }{ \prod _{x\in d(\mu )} (1
- t^{1+l(x)}q^{-a(x)}) (1 - t^{-l(x)}q^{1+a(x)}) }.
\end{equation}
\end{thm}

We briefly review some of the consequences of this formula, as
developed in \cite{GaHa96}.  First, there is reformulation of
\eqref{e:FRn} along the lines of \eqref{e:FR(n,l)-delta}.  Let $\nabla
$ be the linear operator on $\Lambda $ defined by
\begin{equation}\label{e:nabla}
\nabla \Htild _{\mu }(z;q,t) = t^{n(\mu )}q^{n(\mu ')}\Htild _{\mu
}(z;q,t),
\end{equation}
with $n(\mu )$ as in \eqref{e:Htild} and $\mu '$ denoting the
conjugate partition.

\begin{prop}\label{p:FRn-nabla}
The formula \eqref{e:FRn} may be simply expressed as
\begin{equation}\label{e:FRn-nabla}
\Fcal _{R_{n}}(z;q,t) = \nabla e_{n}(z).
\end{equation}
\end{prop}


Next, making use of the known specializations of $\Htild _{\mu
}(z;q,t)$ at $t=q^{-1}$ and $t=1$, we were able to determine the
corresponding specializations of \eqref{e:FRn-nabla}.

\begin{prop}\label{p:t=1/q}
For $t = q^{-1}$ we have
\begin{equation}\label{e:t=1/q}
\begin{aligned}
q^{\binom{n}{2}} \Fcal _{R_{n}}(z;q,q^{-1})
&	 = \frac{1}{1+q+\cdots
+q^{n}}h_{n}\left[Z\smallfrac{1-q^{n+1}}{1-q} \right] \\
&	 = \sum
_{|\lambda | = n} \frac{s_{\lambda }(1,q,\ldots,q^{n})}{1+q+\cdots
+q^{n}}s_{\lambda }(z)
\end{aligned}
\end{equation}
and hence
\begin{equation}\label{e:t=1/q-H}
q^{\binom{n}{2}} \Hcal _{R_{n}}(q,q^{-1}) = (1+q+\cdots +q^{n})^{n-1}.
\end{equation}
In particular, setting $q=1$, we have
\begin{equation}\label{e:dimRn}
\dim R_{n} = (n+1)^{n-1}.
\end{equation}
\end{prop}


The specialization at $t=1$ is most conveniently expressed
combinatorially, in terms of parking functions.  A function $f\colon
\{1,\ldots,n \}\rightarrow \{1,\ldots,n \}$ is called a {\it parking
function} if $|f^{-1}(\{1,\ldots,k \})|\geq k$, for all $1\leq k\leq
n$.  To understand the name, picture a one-way street with $n$ parking
spaces numbered $1$ through $n$.  Suppose that $n$ cars arrive in
succession, each with a preferred parking space given by $f(i)$ for
the $i$-th car.  Each driver proceeds directly to his or her preferred
space and parks there, or in the next available space, if the desired
space is already taken.  The necessary and sufficient condition for
everyone to park without being forced to the end of the street is that
$f$ is a parking function.  The {\it weight} of $f$ is the quantity
$w(f) = \sum _{i=1}^{n}f(i) - i$.  It measures the quantity of
frustration experienced by the drivers in having to pass up occupied
parking spaces.  The symmetric group acts on the set $\PF _{n}$ of
parking functions by permuting the cars (that is, the domain of $f$)
and this action preserves the weight.  Let $\CC \PF _{n} = \bigoplus_{d}
\CC \PF _{n,d}$ be the permutation representation on parking functions,
graded by weight, {\it i.e.}, $\PF _{n,d} = \{f\in \PF _{n}:w(f) = d \}$.

\begin{prop}\label{p:t=1}
For $t = 1$, we have
\begin{equation}\label{e:t=1}
\Fcal R_{n}(z;q,1) = \sum _{d} q^{d}\phi \ch (\varepsilon \otimes \CC
\PF _{n,d}),
\end{equation}
where $\varepsilon $ is the sign representation.  In other words,
$R_{n}$ and $\varepsilon \otimes \CC \PF _{n}$ are isomorphic as singly
graded $S_{n}$-modules when we consider only the $y$-degree in $R_{n}$
and ignore the $x$-degree.
\end{prop}


Since it is known that $|\PF _{n}| = (n+1)^{n-1}$, we again recover the
dimension formula \eqref{e:dimRn}.  Of particular interest is the
subspace $R_{n}^{\epsilon }$ of $S_{n}$-alternating coinvariants,
whose Hilbert series is given by
\begin{equation}\label{e:HRn-alt}
\Hcal _{R_{n}^{\epsilon }}(q,t) = \langle e_{n}(z),\Fcal
_{R_{n}}(z;q,t) \rangle.
\end{equation}
We can expand this by substituting into \eqref{e:FRn} the known
identity
\begin{equation}\label{e:K1n,mu}
\langle e_{n}(z),\Htild _{\mu }(z;q,t) \rangle = \Ktild _{(1^{n}),\mu
}(q,t) = t^{n(\mu )}q^{n(\mu ')},
\end{equation}
obtaining the following result.

\begin{cor}\label{c:Catalan}
The Hilbert series of the $S_{n}$-alternating diagonal coinvariants is
given by 
\begin{equation}\label{e:Cn(qt)}
\Hcal _{R_{n}^{\epsilon }}(q,t) = C_{n}(q,t) \defeq \sum _{|\mu | =
n} \frac{ t^{n(\mu )}q^{n(\mu ')}(1-q)(1-t) \Pi _{\mu }(q,t) B_{\mu
}(q,t) }{ \prod _{x\in d(\mu )} (1 - t^{1+l(x)}q^{-a(x)}) (1 -
t^{-l(x)}q^{1+a(x)}) }.
\end{equation}
\end{cor}


The quantity $C_{n}(q,t)$, studied in \cite{GaHa96, Hai98}, is called
the {\it $q,t$-Catalan polynomial}.  From either
Proposition~\ref{p:t=1/q} or \ref{p:t=1}, we see that $C_{n}(q,t)$ is
a $q,t$-analog of the Catalan number
\begin{equation}\label{e:Cn}
C_{n}(1,1) = \frac{1}{n+1}{\binom{2n}{n}}
\end{equation}
By the corollary above, we have $C_{n}(q,t)\in \NN [q,t]$.
Recently, Garsia and Haglund also proved this by establishing the
following combinatorial interpretation.

\begin{prop}[\cite{GaHag01}]\label{p:GaHag}
Let $D_{n}$ be the set of non-negative integer sequences
$(e_{1},e_{2},\ldots,e_{n})\in \NN ^{n}$ satisfying $e_{1}=0$ and
$e_{k+1}\leq e_{k}+1$ for all $k$.  Put $|e| = \sum _{i}e_{i}$ and let
$i(e)$ be the number of index pairs $i<j$ such that $e_{j} = e_{i}$ or
$e_{j} = e_{i}-1$.  Then
\begin{equation}\label{e:GaHag}
C_{n}(q,t) = \sum _{e\in D_{n}} t^{|e|}q^{i(e)}.
\end{equation}
\end{prop}

We remark that \eqref{e:K1n,mu} has a direct geometric interpretation.
The bundle $P$ is a quotient of $B^{\otimes n}$ (see \cite[Section
3.7]{Hai01}), so we have an equivariant isomorphism of line bundles
\begin{equation}\label{e:Pepsilon}
P_{(1^{n})} = P^{\epsilon }\cong \ext ^{n}B\cong \Ocal (1).
\end{equation}
Hence $\Ktild _{(1^{n}),\mu }(q,t)$, which is the $\TT ^{2}$ character
of the fiber $P_{(1^{n})}(I_{\mu }) = \ext ^{n}B(I_{\mu })$, is equal
to $\prod _{(r,s)\in d(\mu )}t^{r}q^{s} = t^{n(\mu )}q^{n(\mu ')}$.
The notation $\Ocal (1)$ here refers to the very ample line bundle
coming from the projective embedding of $H_{n}$ over $S^{n}\CC ^{2}$
constructed in \cite[Proposition 2.6]{Hai98}.  The identity $\ext
^{n}B\cong \Ocal (1)$ is \cite[Proposition 2.12]{Hai98}.  See
also Proposition~\ref{p:omega-X}, below.


\subsection*{Other character formulas}

The ring $R(n,l)$ and its quotient $R(n,l)/\mfrak R(n,l)$ have $S_{l}$
actions permuting the coordinates $a_{1},b_{1},\ldots,a_{l},b_{l}$,
and commuting with the $S_{n}$ action.  Under our identification of
these rings with the spaces of global sections $H^{0}(H_{n},{P\otimes
B^{\otimes l}})$ and $H^{0}(Z_{n},{P\otimes B^{\otimes l}})$, the
$S_{l}$ action corresponds to permutation of the tensor factors in
$B^{\otimes l}$.

Recall that the {\it Schur functor} $S^{\nu }$ for $\nu $ a partition
of $l$ is defined by 
\begin{equation}\label{e:Schur}
S^{\nu }(W) = (W^{\otimes l})_{\nu } = \Hom ^{S_{n}}(V^{\nu
},W^{\otimes l}).
\end{equation}
It makes sense as a functor on vector spaces and also on vector
bundles.  The following classical result of Schur \cite{Sch27} can be
viewed as a formulation of Schur-Weyl duality.

\begin{prop}\label{p:Schur-Weyl}
If $\alpha \in \operatorname{End}(W)$ has eigenvalues
$t_{1},\ldots,t_{d}$, then the trace of $S^{\nu }(\alpha )\in
\operatorname{End}S^{\nu }(W)$ is given by the Schur function
\begin{equation}\label{e:Schur-Weyl}
s_{\nu }(t_{1},\ldots,t_{d}).
\end{equation}
\end{prop}

\begin{cor}\label{c:SnuB(Imu)}
The Hilbert series of $S^{\nu }(B(I_{\mu }))$ is given by
\begin{equation}\label{e:SnuB(Imu)}
\Hcal _{S^{\nu }(B(I_{\mu }))} = s_{\nu }[B_{\mu }(q,t)]
\end{equation}
in the notation of \eqref{e:f[Y]}.
\end{cor}


Proceeding as in the derivation of Theorems~\ref{t:FR(n,l)} and
\ref{t:FRn}, one obtains the following refinement, which takes account
of the $S_{l}$ action.

\begin{thm}\label{t:(FR)nu}
The Frobenius series of $R(n,l)_{\nu } = \Hom ^{S_{l}}(V^{\nu },
R(n,l))$ is given by
\begin{equation}\label{e:(FR)nu}
\Fcal _{R(n,l)_{\nu }}(z;q,t) = \sum _{|\mu | = n} \frac{ s_{\nu
}[B_{\mu }(q,t)] \Htild _{\mu }(z;q,t) }{ \prod _{x\in d(\mu )} (1 -
t^{1+l(x)}q^{-a(x)}) (1 - t^{-l(x)}q^{1+a(x)}) }.
\end{equation}
Setting $S(n,l,\nu ) = (R(n,l)/\mfrak R(n,l))_{\nu }$, its Frobenius series
 is given by
\begin{equation}\label{e:(FR/m)nu}
\Fcal _{S(n,l,\nu )}(z;q,t) = \sum _{|\mu | = n} \frac{ (1-q)(1-t)
\Pi _{\mu }(q,t) B_{\mu }(q,t) s_{\nu }[B_{\mu }(q,t)] \Htild _{\mu
}(z;q,t) }{ \prod _{x\in d(\mu )} (1 - t^{1+l(x)}q^{-a(x)}) (1 -
t^{-l(x)}q^{1+a(x)}) }.
\end{equation}
\end{thm}

It is convenient to express these identities with the aid of operators
$\Delta _{f}$ defined for any symmetric function $f$ by
\begin{equation}\label{e:nabla-f}
\Delta _{f}\Htild _{\mu }(z;q,t) = f[B_{\mu }(q,t)]\Htild (z;q,t).
\end{equation}
In this notation, the operator $\Delta $ in \eqref{e:Delta} is $\Delta
_{e_{1}}$, and $\nabla $ in \eqref{e:nabla} is the operator which
coincides with $\Delta _{e_{n}}$ in degree $n$, for each $n$.  From
the expressions for the $l=0$ cases of \eqref{e:(FR)nu} and
\eqref{e:(FR/m)nu} in Lemma~\ref{l:R(n,0)} and
Proposition~\ref{p:FRn-nabla}, we get the following corollary.


\begin{cor}\label{c:(FR)nu}
The two Frobenius series in \eqref{e:(FR)nu} and \eqref{e:(FR/m)nu}
may be simply expressed as
\begin{gather}\label{e:(FR)nu-simple}
\Fcal _{R(n,l)_{\nu }}(z;q,t) = \Delta _{s_{\nu }}h_{n}\left[\smallfrac{Z}{(1-q)(1-t)} \right],\\
\label{e:(FR/m)nu-simple}
\Fcal _{S(n,l,\nu )}(z;q,t) = \Delta _{s_{\nu }}\nabla e_{n}(z) =
\Delta _{e_{n}s_{\nu }} e_{n}(z).
\end{gather}
In particular, the expression on the right-hand side is a $q,t$-Schur
positive formal power series in \eqref{e:(FR)nu-simple} and polynomial
in \eqref{e:(FR/m)nu-simple}.
\end{cor}

The operators $\Delta _{s_{\nu }}$ were studied in \cite{BeGaHaTe99},
where we made the following conjecture.

\begin{conj}\label{conj:BGHT+}
The quantity $\Delta _{s_{\nu }}e_{n}(z)$ is a $q,t$-Schur positive
polynomial for all $\nu $ and $n$.
\end{conj}

This statement is stronger than the positivity of the expression in
\eqref{e:(FR/m)nu-simple}, because $\Delta _{f}$ is linear in $f$, and
$e_{n}s_{\nu }$ is a positive linear combination of Schur functions.


\begin{prop}\label{p:BGHT+=chi}
We have
\begin{equation}\label{e:BGHT+=chi}
\chi \Fcal _{\Ocal _{Z_{n}}\otimes P^{*}\otimes S^{\nu }(B)} = \Delta
_{s_{\nu }} e_{n}(z).
\end{equation}
\end{prop}

\begin{proof}
Equation \eqref{e:(FR/m)nu-simple} gives $\chi \Fcal _{\Ocal
_{Z_{n}}\otimes P \otimes S^{\nu }(B)} = \Delta _{s_{\nu }} \nabla
e_{n}(z)$.  To remove the extra factor $\nabla $, we should divide the
numerator in \eqref{e:(FR/m)nu} by $t^{n(\mu )}q^{n(\mu ')}$.  By the
remarks following \eqref{e:Pepsilon}, this is achieved if we replace
$P$ with ${\Ocal (-1)\otimes P}$.  The latter is isomorphic to the
dual bundle $P^{*}$ \cite[eq.~(45)]{Hai01}.
\end{proof}

From this we see that $\Delta _{s_{\nu }}e_{n}(z)$ is at least a
polynomial and that Conjecture~\ref{conj:BGHT+} would be a consequence
of the following strengthening of Theorems~\ref{t:van} and
\ref{t:van-zero}. 

\begin{conj}\label{conj:strongvan}
We have $H^{i}(H_{n},{P^{*}\otimes B^{\otimes l}}) = 0$ for all $i>0$,
and hence also $H^{i}(Z_{n},{P^{*}\otimes B^{\otimes l}}) = 0$ for all
$i>0$.
\end{conj}

Note that the ``hence also'' part follows precisely as in the
derivation of Theorem~\ref{t:van-zero} from Theorem~\ref{t:van}.  The
identification of the spaces of global sections seems rather
difficult, and will not be addressed here.


\section{The operator conjecture}
\label{s:op}

In \cite{Hai94}, we proved the following proposition and conjectured
that the theorem stated below it holds.  The theorem was called the
{\it operator conjecture} in \cite{Hai94}.

\begin{prop}\label{p:op-def}
The space $\Dh _{n}$ of diagonal harmonics defined in \eqref{e:DHn} is
closed under the action of the {\em polarization operators}
\begin{equation}\label{e:Ek}
E_{k} = \sum _{i=1}^{n} y_{i}\partial x_{i}^{k},\quad k>0.
\end{equation}
\end{prop}

\begin{thm}\label{t:op}
The Vandermonde determinant $\Delta (\xx )$ generates $\Dh _{n}$ as a
module for the algebra of operators $\CC [\partial
x_{1},\ldots,\partial x_{n}, E_{1},\ldots,E_{n-1}]$.
\end{thm}


Note that the operators $\partial x_{j}$ and $E_{k}$ all commute, and
that we need not go past $E_{n-1}$, as $E_{k}\Delta (\xx )=0$ for
$k\geq n$.  We will prove the theorem using the isomorphism
\begin{equation}\label{e:psi-1-Rn}
\psi _{1}\colon R_{n}\rightarrow H^{0}(Z_{n},P)
\end{equation}
given by the case $l=0$ of Theorem~\ref{t:van-zero}, where $R_{n}$ is
the ring of diagonal coinvariants.  The first step is to recast
Theorem~\ref{t:op} in ideal-theoretic terms.  There is a symmetric
inner product $(\cdot ,\cdot )$ on $\CC [\xx ,\yy ]$ defined by
\begin{equation}\label{e:apolar}
(f,g) = \left. g(\partial \xx ,\partial \yy )f(\xx ,\yy ) \right|_{\xx
,\yy \mapsto 0}.
\end{equation}
The set of all monomials $\xx ^{h}\yy ^{k}$ is an orthogonal basis,
with $(\xx ^{h}\yy ^{k}, \xx ^{h}\yy ^{k}) = \prod _{i=1}^{n}
(h_{i})!(k_{i})!$.  In particular, this verifies that $(\cdot ,\cdot
)$ is in fact symmetric.  The inner product is compatible with the
grading and non-degenerate.  Since $\CC [\xx ,\yy ]_{d}$ is
finite-dimensional in each degree $d$, we have $I^{\perp \perp } = I$
for any homogeneous subspace $I\subseteq \CC [\xx ,\yy ]$.  One sees
easily from \eqref{e:apolar} that the operator $\partial x_{j}$ is
adjoint to multiplication by $x_{j}$, and likewise for $y_{j}$.  A
polynomial $f$ is orthogonal to an ideal $(g_{1},\ldots,g_{k})$ if and
only if
\begin{equation}\label{e:fperp(g)}
\left. p(\partial \xx ,\partial \yy )g_{i}(\partial \xx ,\partial \yy
) f(\xx ,\yy ) \right|_{\xx ,\yy \mapsto 0} = 0
\end{equation}
for all $i$ and all $p\in \CC [\xx ,\yy ]$.  By Taylor's theorem, this
is equivalent to $g_{i}(\partial \xx ,\partial \yy ) f(\xx ,\yy ) = 0$
for all $i$.  Setting 
\begin{equation}\label{e:anotherI}
I = \mfrak \CC [\xx ,\yy ],
\end{equation}
we therefore see that $\Dh _{n} = I^{\perp }$, or $I = \Dh _{n}^{\perp
}$.  The following version of Theorem~\ref{t:op} in one set of
variables is classical.


\begin{prop}[\cite{Ste64}]\label{p:Steinberg}
Let $I_{0}\subseteq \CC [\xx ]$ be the ideal generated by the
homogeneous maximal ideal in $\CC [\xx ]^{S_{n}}$, or equivalently by
the elementary symmetric functions $e_{1}(\xx ),\ldots,e_{n}(\xx )$,
so that $I_{0}^{\perp }$ is the space of harmonics for the usual
action of $S_{n}$ on $\CC ^{n}$.  Then the Vandermonde determinant
$\Delta (\xx )$ generates $I_{0}^{\perp }$ as a $\CC [\partial \xx
]$-module.
\end{prop}

Returning to the diagonal situation, set
\begin{equation}\label{e:OPn}
\OP _{n} = \CC [\partial \xx , E_{1},\ldots,E_{n-1}]\Delta (\xx ).
\end{equation}
We have $\OP _{n}\subseteq \Dh _{n}$, and hence
\begin{equation}\label{e:IinOP-perp}
I\subseteq \OP _{n}^{\perp },
\end{equation}
and we are to prove that equality holds here.


\begin{prop}\label{p:op-first}
We have $f(\xx ,\yy )\in \OP _{n}^{\perp }$ if and only if
\begin{equation}\label{e:finOPn}
f(\xx ,\phi _{\lambda }(\xx ))\in (e_{1}(\xx ),\ldots,e_{n}(\xx )),
\end{equation}
identically in $\lambda $, where $\phi _{\lambda }(z) = \lambda
_{n-1}z^{n-1} + \cdots + \lambda _{1}z$ is the polynomial of degree
$n-1$ in one variable with zero constant term and generic
coefficients.
\end{prop}

\begin{proof}
Since the adjoint of $\partial x_{j}$ is $x_{j}$, and the adjoint of
$E_{k}$ is $E_{k}^{*} = \sum _{i} x_{i}^{k}\partial y_{i}$, it follows
that we have $f\in \OP _{n}^{\perp }$ if and only if 
\begin{equation}\label{e:op-cond-I}
\Delta (\xx )\perp \CC [\xx ,E_{1}^{*},\ldots,E_{n-1}^{*}]f.
\end{equation}
The formal series $\exp(\lambda _{n-1}E_{n-1}^{*}+\cdots +\lambda
_{1}E_{1}^{*})$ may be viewed as a generating function in the
indeterminates $\lambda _{k}$ for all monomials in the operators
$E_{k}^{*}$.  Condition \eqref{e:op-cond-I} is then equivalent to
\begin{equation}\label{e:op-cond-II}
\exp({\textstyle \sum _{k}\lambda _{k}E_{k}^{*}})f\subseteq
( \CC [\partial \xx ] \Delta (\xx ) )^{\perp }
\end{equation}
holding identically in $\lambda $.  This last condition depends only
on the $y$-degree zero part of $\exp(\sum _{k}\lambda
_{k}E_{k}^{*})f$, so from Proposition~\ref{p:Steinberg} we see that it
is in turn equivalent to
\begin{equation}\label{e:op-cond-III}
(\exp({\textstyle \sum _{k}\lambda _{k}E_{k}^{*}})f)|_{\yy \mapsto
0}\in (e_{1}(\xx ),\ldots,e_{n}(\xx )).
\end{equation}
By Taylor's theorem, $\exp(\lambda _{k}E_{k}^{*})f$ is equal to the
result of substituting $y_{j}+\lambda _{k}x_{j}^{k}$ for $y_{j}$ in
$f$, for all $j$.  Hence $(\exp(\sum _{k}\lambda
_{k}E_{k}^{*})f)|_{\yy \mapsto 0} = f(\xx ,\phi _{\lambda }(\xx ))$,
and the proposition is proved.
\end{proof}


Theorem~\ref{t:op} is a corollary to the preceding proposition and the
next.

\begin{prop}\label{p:op-final} If $f\in \CC [\xx ,\yy ]$ satisfies
$f(\xx ,\phi _{\lambda }(\xx ))\in (e_{1}(\xx ),\ldots,e_{n}(\xx ))$,
with $\phi _{\lambda }$ as in Proposition~\ref{p:op-first}, then
$f(\xx ,\yy )\in \mfrak \CC [\xx ,\yy ]$, where $\mfrak = \CC [\xx
,\yy ]^{S_{n}}_{+}$.
\end{prop}

\begin{proof}
Using Theorem~\ref{t:van-zero}, it suffices to show that the global
section $\psi f(\xx ,\yy )\in H^{0}(H_{n},P)$ restricts to zero on
$Z_{n}$.  Equivalently, we are to show that the function $f(\xx ,\yy
)$ on $X_{n}$ belongs to the ideal of the scheme-theoretic preimage
$\rho ^{-1}(Z_{n})$.


Let $U_{x}\subseteq H_{n}$ be the open set consisting of ideals $I$
such that $x$ generates the tautological fiber $B(I) = \CC [x,y]/I$ as
a $\CC $-algebra, that is,
\begin{equation}\label{e:Ux}
U_{x} = \{I\in H_{n}: \text{$\{1,x,\ldots,x^{n-1} \}$ is a basis of
$B(I)$} \}.
\end{equation}
As shown in \cite[Section~3.6]{Hai01}, $U_{x}$ is an affine cell with
coordinates $e_{1},\ldots,e_{n}$, $\gamma _{0},\ldots,\gamma _{n-1}$
such that the equations of the universal family over $U_{x}$ are given
in terms of these and the coordinates $x,y$ on $\CC ^{2}$ by
\begin{equation}\label{e:F-over-Ux}
\begin{gathered}
x^{n} - e_{1}x^{n-1}+\cdots + (-1)^{n}e_{n} = 0\\
y = \gamma _{n-1}x^{n-1}+\cdots +\gamma _{1}x+\gamma _{0}.
\end{gathered}
\end{equation}
The preimage $U'_{x} = \rho ^{-1}(U_{x})$ of $U_{x}$ in $X_{n}$ is an
affine cell with coordinates $x_{1},\ldots,x_{n}, \gamma
_{0},\ldots,\gamma _{n}$.  The morphism $\rho \colon X_{n}\rightarrow
H_{n}$ is given on the coordinate level by the identification of
$e_{i}$ with the $i$-th elementary symmetric function $e_{i}(\xx )$.
Each coordinate pair $x_{j}, y_{j}$ on $X_{n}$ satisfies equations
\eqref{e:F-over-Ux}, so the coordinates $y_{j}$ are given in terms of
$\xx ,\ggam $ by $y_{j} = \phi _{\gamma }(x_{j})$, where $\phi
_{\gamma }(z)$ is the polynomial $\gamma _{n-1}z^{n-1} + \cdots +
\gamma _{1}z +\gamma _{0}$ with coefficients $\ggam $.

The zero fiber $Z_{n}$ is irreducible \cite{Bri77}, so $U_{x}\cap
Z_{n}$ is dense in $Z_{n}$, and it suffices to check that the section
represented by $f$ is zero there.  In terms of the coordinates $\ee
,\ggam $ on $U_{x}$, the ideal of $U_{x}\cap Z_{n}$ is
$(\gamma _{0},\ee )$, so the coordinate ring of the scheme-theoretic
preimage  $U_{x}'\cap \rho ^{-1}(Z_{n})$ is 
\begin{equation}\label{e:rho-inv-Zn}
\CC [\xx ,\gamma _{1},\ldots,\gamma _{n-1}]/(e_{1}(\xx
),\ldots,e_{n}(\xx )).
\end{equation}
In terms of the coordinates $\xx ,\ggam $, the given function $f(\xx
,\yy )$ becomes $f(\xx ,\phi _{\gamma }(\xx ))$, which belongs to
$(e_{1}(\xx ),\ldots,e_{n}(\xx ))$ by hypothesis.
\end{proof}


\section{Proof of the main theorem}
\label{s:van}

We will prove Theorem~\ref{t:van} by combining two results
from \cite{Hai01}---the isomorphism $\CC ^{2n}\sslash S_{n}\cong
H_{n}$ and the theorem that $R(n,l)$ is a free $\CC [\yy
]$-module---with the theorem of Bridgeland, King and Reid mentioned in
the introduction.  We begin by reviewing these results.

Let $V = \CC ^{m}$ be a complex vector space and $G$ a finite subgroup
of $\SL (V)$.  As in Section~\ref{s:main}, we have a diagram
\begin{equation}\label{e:M//G-diagram}
\begin{CD}
X&		@>{f}>>&	V\\
@V{\rho }VV&	&		@VVV\\
V\sslash G&	@>{\sigma }>>&	V/G,
\end{CD}
\end{equation}
whose special case for $V=\CC ^{2n}$, $G = S_{n}$ is
\eqref{e:X-H-diagram}.  Let $D(V\sslash G)$ be the derived category of
complexes of sheaves of $\Ocal _{V\sslash G}$-modules with bounded,
coherent cohomology, and $D^{G}(V)$ the derived category of complexes
of $G$-equivariant sheaves of $\Ocal _{V}$-modules, again with
bounded, coherent cohomology.  Bridgeland, King and Reid define a
functor
\begin{equation}\label{e:Phi}
\Phi \colon D(V\sslash G)\rightarrow D^{G}(V)
\end{equation}
by the formula
\begin{equation}\label{e:Phi-defn}
\Phi  = Rf_{*}\circ \rho ^{*}.
\end{equation}
Note that $\rho $ is flat, so we can write $\rho ^{*}$ instead of
$L\rho ^{*}$ here.


\begin{thm}[\cite{BrKiRe01}]\label{t:BKR}
Suppose that the Chow morphism $V\sslash G\rightarrow V/G$ satisfies
the following smallness criterion: for every $d$, the locus of points
$x\in V/G$ such that $\dim \sigma^{-1}(x)\geq d$ has codimension at
least $2d-1$.  Then 
\begin{itemize}
\item [(1)] $V\sslash G$ is a crepant resolution of singularities of $V/G$,
{\it i.e.}, it is non-singular and its canonical line bundle is
trivial, and
\item [(2)] the functor $\Phi $ is an equivalence of categories.
\end{itemize}
\end{thm}

We apply the theorem with $V = \CC ^{2n}$ and $G = S_{n}$.  Note that
$S_{n}$, acting diagonally, is a subgroup of $\SL(\CC ^{2n})$.  It is
known \cite{Hai01,Nak99} that $\omega _{H_{n}} \cong \Ocal _{H_{n}}$,
so $\CC ^{2n}\sslash S_{n}\cong H_{n}$ is a crepant resolution of $\CC
^{2n}/S_{n} = S^{n}\CC ^{2}$.  Moreover, the smallness criterion in
Theorem~\ref{t:BKR} holds.  This follows either from the description
of the fibers of the Chow morphism due to Brian\c{c}on \cite{Bri77},
or from the observation in \cite{BrKiRe01} that, conversely to
Theorem~\ref{t:BKR}, the criterion holds whenever $G$ preserves a
symplectic form on $V$ and $V\sslash G$ is a crepant resolution.  We
identify $D^{S_{n}}(\CC ^{2n})$ with the derived category of bounded
complexes of finitely-generated $S_{n}$-equivariant $\CC [\xx ,\yy
]$-modules.  The functor $Rf_{*}$ is thereby identified with $R\Gamma
_{X_{n}}$.  Since $\rho $ is finite and therefore affine, and $P =
\rho _{*}\Ocal _{X_{n}}$, the functor $R\Gamma _{X_{n}}\circ \rho
^{*}$ is naturally isomorphic to $R\Gamma _{H_{n}}({P\otimes -})$.

\begin{cor}\label{c:BKRv2}
The functor $\Phi = R\Gamma ({P\otimes -})$ is an equivalence of
categories $\Phi \colon D(H_{n})\rightarrow D^{S_{n}}(\CC ^{2n})$.
\end{cor}


Using this we can reformulate our main theorem.

\begin{prop}\label{p:van-via-Phi}
Theorem~\ref{t:van} is equivalent to the identity in $D^{S_{n}}(\CC
^{2n})$
\begin{equation}\label{e:Phi=R(n,l)}
\Phi B^{\otimes l} \cong  R(n,l),
\end{equation}
where the isomorphism is given by the map $R(n,l)\rightarrow \Phi
B^{\otimes l}$ obtained by composing the canonical natural
transformation $\Gamma \rightarrow R\Gamma $ with the homomorphism
$\psi $ in \eqref{e:R(n,l)->H^0}.
\end{prop}

We will prove identity \eqref{e:Phi=R(n,l)}, and thus
Theorem~\ref{t:van}, by using the {\it inverse} Bridgeland--King--Reid
functor $\Psi \colon D^{S_{n}}(\CC ^{2n})\rightarrow D(H_{n})$, which
also has a simple description in our case.  In general, as observed
in \cite{BrKiRe01}, the inverse functor $\Psi $ can be calculated
using Grothendieck duality as the right adjoint of $\Phi $, given by
the formula
\begin{equation}\label{e:Psi-general}
\Psi  = (\rho _{*}(\omega _{X}\overset{L}{\otimes }Lf^{*}-))^{G}.
\end{equation}
To simplify this, we use the following result from \cite{Hai01}.

\begin{prop}\label{p:omega-X}
The line bundle $\Ocal (1) = \ext ^{n}B$ is the Serre twisting sheaf
induced by a natural embedding of $H_{n}$ as a scheme projective over
$S^{n}\CC ^{2}$.  Writing $\Ocal (1)$ also for its pullback to
$X_{n}$, we have that $X_{n}$ is Gorenstein with canonical sheaf
$\omega _{X_{n}} \cong \Ocal (-1)$.
\end{prop}


We need an extra bit of information not contained in the proposition.
There are two possible equivariant $S_{n}$ actions on $\Ocal
_{X_{n}}(1)$: the trivial action coming from the definition of $\Ocal
_{X_{n}}(1)$ as $\rho ^{*}\Ocal _{H_{n}}(1)$, or its twist by the sign
character of $S_{n}$.  The latter action is the correct one, in the
sense that the isomorphism $\omega _{X_{n}}\cong \Ocal (-1)$ is
$S_{n}$-equivariant for this action, as can be seen from the proof
in \cite{Hai01}.  Taking this into account, and using the fact that
$\Ocal _{X_{n}}(-1)$ is pulled back from $H_{n}$, we have the
following description of the inverse functor.

\begin{prop}\label{p:Psi}
The inverse of the functor $\Phi $ in Corollary \ref{c:BKRv2} is given
by 
\begin{equation}\label{e:Psi-X}
\Psi = \Ocal (-1)\otimes (\rho _{*}\circ  Lf^{*})^{\epsilon }.
\end{equation}
Here $(-)^{\epsilon }$ denotes the functor of $S_{n}$-alternants, {\it
i.e.}, $A^{\epsilon } = \Hom ^{S_{n}}(\varepsilon ,A)$, where
$\varepsilon $ is the sign representation.
\end{prop}


Now we recall the algebraic result that was the key technical tool
in \cite{Hai01}.

\begin{thm}\label{t:polygraph}
The polygraph coordinate ring $R(n,l)$ is a free $\CC [\yy ]$-module.
\end{thm}

We need to strengthen this in two ways.  Any automorphism of $\CC
^{2}$ induces an automorphism of $\CC ^{2n+2l}$, and the corresponding
automorphism of $\CC [\xx ,\yy ,\aa,\bb ]$ leaves invariant the
defining ideal $I(n,l)$ of $Z(n,l)$.  In particular, this is so for
translations in the $x$-direction, which also leave invariant the
ideal $(\yy )$ and hence $I(n,l)+(\yy )$.  This implies that any of
the coordinates $x_{i}$, $a_{i}$ is a non-zero-divisor in $R(n,l)/(\yy
)$, yielding the following two corollaries.

\begin{cor}\label{c:polygraph}
The coordinate ring $R(n,l)$ is a free $\CC [x_{1},\yy ]$-module.
\end{cor}

\begin{cor}\label{c:R(n,l)-free-res}
The coordinate ring $R(n,l)$ has a free resolution of length $n-1$ as
a $\CC [\xx ,\yy ]$-module.
\end{cor}


As in \cite[Definition 4.1.1]{Hai01}, the polygraph $Z(n,l)$ is the
union of linear subspaces $W_{f}\subseteq \CC ^{2n+2l}$ defined by
\begin{equation}\label{e:Wf}
W_{f} = V(I_{f}),\quad I_{f} = (a_{i}-x_{f(i)},b_{i}-y_{f(i)}:1\leq
i\leq l)
\end{equation}
for all functions $f\colon \{1,\ldots,l \}\rightarrow \{1,\ldots,n
\}$.  The polygraph ring can be defined with any ground ring $S$ in
place of $\CC $ as
\begin{equation}\label{e:R(n,l)-over-S}
R(n,l) = S[\xx ,\yy ,\aa ,\bb ]/I(n,l),\quad I(n,l) = \bigcap _{f} I_{f},
\end{equation}
with $I_{f}$ as above.  Theorem~\ref{t:polygraph} holds in this more
general setting \cite[Theorem 4.3]{Hai01}.  If $\theta $ is an
automorphism of $S[x,y]$ as an $S$-algebra, then the automorphism
$\theta ^{\otimes (n+l)}$ of $S[\xx ,\yy ,\aa ,\bb ] \cong
S[x,y]^{\otimes (n+l)}$ leaves $I(n,l)$ invariant, inducing an
automorphism of $R(n,l)$.  Hence we have the following corollary.

\begin{cor}\label{c:R(n,l)-free-gamma}
Let $S$ be a $\CC $-algebra and let $y'$ denote the image of $y$ under
some automorphism of $S[x,y]$ as an $S$-algebra.  Then ${S\otimes
_{\CC }R(n,l)}$ is a free $S[y_{1}',\ldots,y_{n}']$-module.
\end{cor}

In addition to the results on polygraphs we need the following local
structure theorem for $X_{n}$.  It allows us to assume by induction on
$n$ that a desired geometric result holds locally over the open locus
consisting of points $I\in H_{n}$ such that $V(I)$ is not concentrated
at a single point of $\CC ^{2}$.


\begin{prop}\label{p:product}
Let $U_{k}\subseteq X_{n}$ be the open set consisting of points
$(I,P_{1},\ldots,P_{n})$ for which $\{P_{1},\ldots,P_{k} \}$ and
$\{P_{k+1},\ldots,P_{n} \}$ are disjoint.  Then $U_{k}$ is isomorphic
to an open set in $X_{k}\times X_{n-k}$.  More precisely, the morphism
$f\colon X_{n}\rightarrow \CC ^{2n}$ restricted to $U_{k}$ corresponds
to the restriction of $f_{k}\times f_{n-k}\colon X_{k}\times
X_{n-k}\rightarrow \CC ^{2k}\times \CC ^{2(n-k)} = \CC ^{2n}$.

The pullback $F_{n}' = F_{n}\times X_{n}\,/\,H_{n}$ of the universal
family to $X_{n}$ decomposes over $U_{k}$ as the disjoint union
$F_{n}' = F'_{k}\times X_{n-k}\cup X_{k}\times F'_{n-k}$ of the
pullbacks of the universal families from $H_{k}$ and $H_{n-k}$.  Hence
the tautological sheaf $\rho ^{*}B$ decomposes as $\rho ^{*}B = \eta
_{k}^{*} \rho _{k}^{*}B_{k}\oplus \eta _{n-k}^{*}\rho
_{n-k}^{*}B_{n-k}$, where $\eta _{k}$ and $\eta _{n-k}$ are the
projections of $X_{k}\times X_{n-k}$ on the factors.
\end{prop}

The final piece of our puzzle will be supplied by a fundamental result
of commutative algebra known as the new intersection theorem.

\begin{thm}[\cite{PeSz74,Rob87,Rob89}]\label{t:NIT}
Let $0\rightarrow C_{n}\rightarrow \cdots \rightarrow C_{1}\rightarrow
C_{0}\rightarrow 0$ be a bounded complex of locally free coherent
sheaves on a Noetherian scheme $X$.  Denote by $\Supp (C\subdot )$ the
union of the supports of the homology sheaves $H_{i}(C\subdot )$.
Then every component of $\Supp (C\subdot )$ has codimension at most
$n$ in $X$.  In particular, if $C\subdot $ is exact on an open set
$U\subseteq X$ whose complement has codimension exceeding $n$, then
$C\subdot $ is exact.
\end{thm}


\noindent {\it Proof of Theorem~\ref{t:van}.}  By
Proposition~\ref{p:van-via-Phi}, we have a map
\begin{equation}\label{e:R(n,l)->PhiBl}
R(n,l)\rightarrow \Phi B^{\otimes l}
\end{equation}
in the derived category $D^{S_{n}}(\CC ^{2n})$, and it suffices to
show that it is an isomorphism.  Applying the inverse functor $\Psi $
yields a map
\begin{equation}\label{e:PsiR(n,l)->Bl}
\Psi R(n,l)\rightarrow B^{\otimes l}
\end{equation}
in $D(H_{n})$, and we can equally well show that this is an
isomorphism.  Let $C$ be the third vertex of a distinguished
triangle
\begin{equation}\label{e:triangle}
C[-1]\rightarrow \Psi R(n,l)\rightarrow B^{\otimes l}\rightarrow C.
\end{equation}
We are to show that $C=0$.

We may compute $\Psi R(n,l)$ as follows.  By Corollary
\ref{c:R(n,l)-free-res}, the $\CC [\xx ,\yy ]$-algebra $R(n,l)$ has a
free resolution of length $n-1$.  We can assume that the resolution is
$S_{n}$-equivariant, for instance by taking a graded minimal free
resolution.  In derived category terminology, we have an
$S_{n}$-equivariant complex of free $\CC [\xx ,\yy ]$-modules
\begin{equation}\label{e:R(n,l)-res}
A\subdot  =\cdots \rightarrow 0\rightarrow A_{n-1}\rightarrow \cdots
\rightarrow A_{1}\rightarrow A_{0}\rightarrow 0\rightarrow \cdots
\end{equation}
quasi-isomorphic to $R(n,l)$.  Using the formula for $\Psi $ from
Proposition \ref{p:Psi}, we have $\Psi R(n,l) = {\Ocal (-1)\otimes
(\rho _{*}f^{*}A\subdot )^{\epsilon }}$.  Moreover, since $\rho $ is
flat, and since the functor $(-)^{\epsilon }$ is a direct summand of
the identity functor, ${\Ocal (-1)\otimes (\rho _{*}f^{*}A\subdot
)^{\epsilon }}$ is a complex of locally free sheaves.  Since
$B^{\otimes l}$ is a sheaf, the map $\Psi R(n,l)\rightarrow B^{\otimes
l}$ in \eqref{e:PsiR(n,l)->Bl} is represented by an honest
homomorphism of complexes, and not merely by a quasi-isomorphism.  The
object $C$ is represented by the mapping cone of this homomorphism,
namely, the complex of locally free sheaves
\begin{equation}\label{e:mapping-cone-C}
0\rightarrow C_{n}\rightarrow \cdots \rightarrow C_{2}\rightarrow
C_{1}\rightarrow B^{\otimes l}\rightarrow 0,
\end{equation}
where $C_{i} = {\Ocal (-1)\otimes (\rho _{*}f^{*}A_{i-1})^{\epsilon
}}$.  We are to prove that this complex is exact.


Let $U\subseteq H_{n}$ be the open set of points $I$ such that $V(I)$
contains at least two distinct points of $\CC ^{2}$.  Let $U_{k}$ be
the open subset in $X_{n}$ on which $\{P_{1},\ldots,P_{k} \}$ is
disjoint from $\{P_{k+1},\ldots,P_{n} \}$.  Clearly the open set $\rho
^{-1}(U)\subseteq X_{n}$ is the union of open sets conjugate by some
permutation $w\in S_{n}$ to $U_{k}$ for some $0<k<n$.  On $U_{k}$, the
decomposition of the tautological sheaf $\rho ^{*}B$ from Proposition
\ref{p:product} induces a decomposition of $\rho ^{*}B^{\otimes l}$ as
a a direct sum
\begin{equation}\label{e:Bl-decomp}
\rho ^{*}B^{\otimes l}\cong \bigoplus _{j=0}^{l}\binom{l}{j}\cdot (\eta
_{k}^{*}\rho _{k}^{*}B_{k})^{\otimes j}\otimes (\eta _{n-k}^{*}\rho
_{n-k}^{*}B_{n-k})^{\otimes l-j}.
\end{equation}

Let $R(n,l)^{\sim }$ be the sheaf of $\Ocal _{\CC ^{2n}}$ modules
corresponding to the $\CC [\xx ,\yy ]$-module $R(n,l)$.  We partition
the set $\{1,\ldots,n \}$ into two subsets $S_{1} = \{1,\ldots,k \}$
and $S_{2} = \{k+1,\ldots,n \}$, and define $\alpha \colon
\{1,\ldots,n \}\rightarrow \{1,2 \}$ to be the function mapping the
elements of $S_{i}$ to $i$.  Let $U_{k}'$ be the open subset
consisting of points $(P_{1},\ldots,P_{n})\in \CC ^{2n}$ satisfying
the same condition that defines $U_{k}$, namely that $\{
P_{1},\ldots,P_{k}\}$ and $\{P_{k+1},\ldots,P_{n} \}$ are disjoint.
Over $U_{k}'$, components $W_{f}$, $W_{g}$ of the polygraph $Z(n,l)$
are disjoint if $\alpha \circ f\not =\alpha \circ g$.  Hence $Z(n,l)$
decomposes over $U_{k}'$ as a union of $2^{l}$ disjoint closed
subschemes $Z_{h}$, indexed by functions $h:\{1,\ldots,l \}\rightarrow
\{1,2 \}$, where $Z_{h}$ is the union of the components $W_{f}$ for
which $\alpha \circ f = h$.  Each subscheme $Z_{h}$ is isomorphic over
$U_{k}'$ to $Z(k,j)\times Z(n-k,l-j)$, where $j = |h^{-1}(\{1 \})|$.
The number of $Z_{h}$ that occur for a given value of $j$ is
$\binom{l}{j}$.  This decomposition of $Z(n,l)$ gives a direct sum
decomposition of $R(n,l)^{\sim }$ on $U_{k}'$ as
\begin{equation}\label{e:R(n,l)-decomp}
R(n,l)^{\sim }\cong \bigoplus _{j=0}^{l}\binom{l}{j}\cdot
R(k,j)^{\sim }\otimes R(n-k,l-j)^{\sim }.
\end{equation}
The decompositions \eqref{e:Bl-decomp} and \eqref{e:R(n,l)-decomp} are
compatible with the map $\psi \colon R(n,l)\rightarrow
H^{0}(X_{n},\rho ^{*}B^{\otimes l})$ in \eqref{e:R(n,l)->H^0}.  More
precisely, they are compatible with the restriction to $U_{k}'$ of the
induced sheaf homomorphism $\psi ^{\sim }\colon R(n,l)^{\sim
}\rightarrow f_{*}\rho ^{*}B^{\otimes l}$.


Now assume by induction that Theorem~\ref{t:van} holds for smaller
values of $n$, the base case $n=1$ being trivial.  The preceding
remarks then show that the map $R(n,l)\rightarrow \Phi B^{\otimes l}$
in \eqref{e:R(n,l)->PhiBl} restricts to an isomorphism on the open set
$U'\subseteq \CC ^{2n}$ of points $(P_{1},\ldots,P_{n})$ with
$P_{1},\ldots,P_{n}$ not all equal.  The functors $\Phi $ and $\Psi $
are defined locally with respect to $S^{n}\CC ^{2}$, so we conclude
that the map $\Psi R(n,l)\rightarrow B^{\otimes l}$ in
\eqref{e:PsiR(n,l)->Bl} is an isomorphism on $U$, and hence the
complex $C$ in \eqref{e:mapping-cone-C} is exact on $U$.  The
complement of $U$ in $H_{n}$ is isomorphic to $\CC ^{2}\times Z_{n}$,
so it has dimension $n+1$ and codimension $n-1$.  Before applying
Theorem~\ref{t:NIT}, we first need to enlarge $U$ to an open set whose
complement has codimension $n+1$.  The desired open set will be $U\cup
U_{x}\cup U_{y}$, where $U_{x}$ is as in \eqref{e:Ux}, and $U_{y}$ is
defined in the obvious analogous way.  Its complement is isomorphic to
$\CC ^{2}\times (Z_{n}\setminus (U_{x}\cup U_{y}))$, which has
codimension $n+1$ by the following lemma.


\begin{lem}\label{l:UxUy}
The complement $Z_{n}\setminus (U_{x}\cup U_{y})$ of $U_{x}\cup U_{y}$
in the zero fiber has dimension $n-3$.
\end{lem}

\begin{proof}
Let $V = Z_{n}\setminus (U_{x}\cup U_{y})$.  Interpreting $\dim V < 0$
to mean that $V$ is empty, the lemma holds trivially for $n=1$, so we
can assume $n\geq 2$.  We consider the decomposition of $Z_{n}$ into
affine cells as in \cite{Bri77, ElSt87}, and show that each cell
intersects $V$ in a locus of dimension at most $n-3$.  There is one
open cell, of dimension $n-1$.  This cell is actually $U_{x}\cap
Z_{n}$, so it is disjoint from $V$.  There is also one cell of
dimension $n-2$.  It has non-empty intersection with $U_{y}$, so its
intersection with $V$ has dimension at most $n-3$.  In fact this
intersection has dimension exactly $n-3$, since the complement of
$U_{y}$ is the zero locus of a section of the line bundle $\ext ^{n}B
= \Ocal (1)$.  All remaining cells have dimension less than or equal
to $n-3$.
\end{proof}

We digress briefly to point out the geometric meaning of this lemma.
For $I$ in the zero fiber, the fiber $B(I)$ is an Artin local $\CC
$-algebra with maximal ideal $(x,y)$.  The point $I$ belongs to
$U_{x}\cup U_{y}$ if and only if the maximal ideal is principal, that
is, $B(I)$ has embedding dimension one, or equivalently, the
corresponding closed subscheme $V(I)$ is a subscheme of some smooth
curve through the origin in $\CC ^{2}$.  In this case $I$ is said to
be {\it curvilinear}.  The lemma says that the non-curvilinear locus
has codimension two in the zero fiber.


The proof of Theorem~\ref{t:van} is now completed by the following
lemma and its symmetric partner with $U_{y}$ in place of $U_{x}$.

\begin{lem}\label{l:Ux-calculation}
The map $\Psi R(n,l)\rightarrow B^{\otimes l}$ restricts to an
isomorphism on $U_{x}$.
\end{lem}

\begin{proof}
Recall the description in the proof of Proposition~\ref{p:op-final} of
the coordinates on $U_{x}$ and its preimage $U_{x}' = \rho
^{-1}(U_{x})$ in $X_{n}$.  The coordinates on $U_{x}'$ are $\xx ,\ggam
$, with $y_{j}$ equal to $\phi _{\gamma }(x_{j})$, where $\phi
_{\gamma }(z) = \gamma _{n-1}z^{n-1} + \cdots + \gamma _{1}z +\gamma
_{0}$.  The coordinates on $U_{x}$ are $\ee ,\ggam $, where $e_{i} =
e_{i}(\xx )$ is the $i$-th elementary symmetric function, so $\CC [\ee
,\ggam ] = \CC [\xx ,\ggam ]^{S_{n}}$.

We have  a trivial isomorphism
\begin{equation}\label{e:C[x,a]=C[x,y,a]/I}
\CC [\xx ,\ggam ]\cong \CC [\xx ,\yy ,\ggam ]
/(y_{j}-\phi _{\gamma }(x_{j})\colon 1\leq j\leq n)
\end{equation}
which is nonetheless useful because it describes $\CC [\xx ,\ggam ]$
as a $\CC [\xx ,\yy ]$-module.  Since $\CC [\xx ,\ggam ]$ and $\CC
[\xx ,\yy ,\ggam ]$ are polynomial rings of dimension $2n$ and $3n$,
respectively, the ideal in \eqref{e:C[x,a]=C[x,y,a]/I} is a complete
intersection ideal.  Hence the Koszul complex $K\subdot (\yy -\phi
_{\gamma }(\xx ))$ over $\CC [\xx ,\yy ,\ggam ]$ is a free resolution
of $\CC [\xx ,\ggam ]$ as a $\CC [\xx ,\yy ,\ggam ]$-module, and
therefore as a $\CC [\xx ,\yy ]$-module.


The restriction of $Lf^{*}R(n,l)$ to the affine open set $U'_{x}$ is
the complex of sheaves associated to the complex of modules $\CC [\xx
,\ggam ]\overset{L}{\underset{\CC [\xx ,\yy ]}{\otimes }} R(n,l)$.
This can be computed by tensoring $R(n,l)$ with the above free
resolution of $\CC [\xx ,\ggam ]$, and the result is the Koszul
complex $K\subdot (\yy -\phi _{\gamma }(\xx ))$ over ${\CC [\ggam
]\otimes _{\CC }R(n,l)}$.  It follows from Corollary
\ref{c:R(n,l)-free-gamma} that this Koszul complex is a free
resolution of ${\CC [\ggam ]\otimes _{\CC }R(n,l)/(\yy -\phi _{\gamma
}(\xx ))}$.  In other words, on $U'_{x}$ we have $Lf^{*}R(n,l) =
f^{*}R(n,l)$, and we have a description of this object as the sheaf
associated to the $\CC [\xx ,\ggam ]$-algebra ${\CC [\ggam ]\otimes
_{\CC }R(n,l)/(\yy -\phi _{\gamma }(\xx ))}$.  It follows that $\Psi
R(n,l) = {\Ocal (-1)\otimes (\rho _{*}Lf^{*}R(n,l))^{\epsilon }}$ is
described on $U_{x}$ as the sheaf associated to the
$S_{n}$-alternating part of this algebra, regarded as a module over
$\CC [\ee ,\ggam ] = \CC [\xx ,\ggam ]^{S_{n}}$.

The equations \eqref{e:F-over-Ux} of the universal family give us the
description of the tautological bundle $B$ as a sheaf of algebras on
$U_{x}$, from which we can get a description of $B^{\otimes l}$.  To
make the variable names match the ones in $R(n,l)$, we should replace
$x,y$ with variables $a_{i},b_{i}$ standing for the generators of the
$i$-th tensor factor in $B^{\otimes l}$.  In this notation,
$B^{\otimes l}$ is the sheaf associated to the $\CC [\ee ,\ggam
]$-algebra
\begin{equation}\label{e:Bl-on-Ux}
\CC [\ee ,\ggam ,\aa ,\bb ]/\sum _{i=1}^{l}(b_{i}-\phi _{\gamma
}(a_{i}), \prod _{j=1}^{n}(a_{i}-x_{j})).
\end{equation}
Note that the products $\prod _{j=1}^{n}(a_{i}-x_{j})$ written here
really only depend on $\aa $ and the elementary symmetric functions
$e_{i} = e_{i}(\xx )$.


The map $\Psi R(n,l)\rightarrow B^{\otimes l}$ is now expressed in
local coordinates as a homomorphism from 
\begin{equation}\label{e:PsiR(n,l)-described}
(\CC [\ggam ]\otimes _{\CC }R(n,l)/(\yy -\phi _{\gamma }(\xx
)))^{\epsilon }
\end{equation}
to the algebra in \eqref{e:Bl-on-Ux}.  The algebra ${\CC [\ggam
]\otimes _{\CC }R(n,l)/(\yy -\phi _{\gamma }(\xx ))}$ is generated by
the variables $\xx ,\ggam ,\aa ,\bb $, all of which are
$S_{n}$-invariant except $\xx $.  It follows that all its
$S_{n}$-alternating elements are multiples of the Vandermonde
determinant $\Delta (\xx )$ by polynomials in $\ggam ,\aa ,\bb $ and
the elementary symmetric functions $e_{i}(\xx )$.  Written out
explicitly, the homomorphism in question sends an element $\Delta (\xx
)p(\ee ,\ggam ,\aa ,\bb )$ to $p(\ee ,\ggam ,\aa ,\bb )$.  The
division by $\Delta (\xx )$ here reflects the presence of the factor
$\Ocal (-1)$ in the formula for $\Psi R(n,l)$.  The space of global
sections of $\Ocal (1)$ on $X_{n}$ can be identified with the ideal
$J\subseteq \CC [\xx ,\yy ]$ generated by $\CC [\xx ,\yy ]^{\epsilon
}$, in such a way that $\Delta (\xx )$ represents the essentially
unique section which vanishes nowhere on $U_{x}$.

Let us denote the above-described homomorphism by $\xi $.  Since the
algebra in \eqref{e:Bl-on-Ux} is generated by the variables $\ee
,\ggam ,\aa ,\bb $, it is clear that $\xi $ is surjective.  The
injectivity of $\xi $ amounts to saying that the expressions
$b_{i}-\phi _{\gamma }(a_{i})$ and $\prod _{j=1}^{n}(a_{i}-x_{j})$ are
annihilated by $\Delta (\xx )$ in ${\CC [\ggam ]\otimes _{\CC
}R(n,l)/(\yy -\phi _{\gamma }(\xx ))}$.  This condition is clearly
necessary, and it is sufficient since it makes multiplication by
$\Delta (\xx )$ a well-defined left inverse to $\xi $.  The products
$\prod _{j=1}^{n}(a_{i}-x_{j})$ are zero in $R(n,l)$ and thus present
no difficulty.  The expressions $b_{i}-\phi _{\gamma }(a_{i})$ are
more subtle, as they do not vanish in ${\CC [\ggam ]\otimes _{\CC
}R(n,l)/(\yy -\phi _{\gamma }(\xx ))}$.


For $x_{1},\ldots,x_{n}$ distinct, the Lagrange
interpolation problem
\begin{equation}\label{e:Lagrange}
y_{j} = \sum _{k=0}^{n-1}\beta _{k}x_{j}^{k},\quad 1\leq j\leq n,
\end{equation}
is solved by a formula giving the coefficients $\beta _{k}$ as
rational functions of the form $\beta _{k} = \Delta _{k}(\xx ,\yy
)/\Delta (\xx )$, where $\Delta _{k}$ is a certain determinant
involving the variables $\xx $, $\yy $.  Multiplying through by
$\Delta (\xx )$ yields the identity of polynomials
\begin{equation}\label{e:yi-Delta}
y_{j}\Delta (\xx ) = \sum _{k=0}^{n-1}\Delta _{k}(\xx ,\yy
) x_{j}^{k},\quad 1\leq j\leq n.
\end{equation}
On each component $W_{f}$ of $Z(n,l)$ we have $a_{i} = x_{f(i)}$,
$b_{i} = y_{f(i)}$, and therefore
\begin{equation}\label{e:bi-Delta}
b_{i}\Delta (\xx ) = \sum _{k=0}^{n-1}\Delta _{k}(\xx ,\yy )
a_{i}^{k},\quad 1\leq i\leq l.
\end{equation}
Since these equations hold on every component of $Z(n,l)$, they hold
identically in $R(n,l)$.  Similarly, for arbitrary values of the
parameters $\ggam $, we may substitute $\phi _{\gamma }(\xx )$ for
$\yy $ in \eqref{e:yi-Delta} and then let $x_{j} = a_{i}$, to obtain
the identity
\begin{equation}\label{e:phi_gamma(a)}
\phi _{\gamma }(a_{i}) \Delta (\xx ) = \sum _{k=0}^{n-1}\Delta
_{k}(\xx ,\phi _{\gamma }(\xx ))a_{i}^{k},\quad 1\leq i\leq l,
\end{equation}
valid when $a_{i}$ is equal to any of the $x_{j}$.  Again this holds
on every component of $Z(n,l)$ and hence as an identity in $R(n,l)$.
Subtracting \eqref{e:phi_gamma(a)} from \eqref{e:bi-Delta}, we see
that $\Delta (\xx )$ annihilates $b_{i}-\phi _{\gamma }(a_{i})$ in
${\CC [\ggam ]\otimes _{\CC }R(n,l)/(\yy -\phi _{\gamma }(\xx ))}$,
which was the only thing left to prove.
\end{proof}


\bibliographystyle{hamsplain}
\bibliography{../references}

\end{document}